\documentclass{amsart}
\usepackage{amsthm,amsmath,amssymb,amscd,graphics,psfig,epic,eepic}
{\end{list}}
{
   \newtheorem{theorem}{Theorem}[subsection]                     
   \newtheorem{proposition}[theorem]{Proposition}     
   \newtheorem{lemma}[theorem]{Lemma}

   \newtheorem{corollary}[theorem]{Corollary}
   \newtheorem{conjecture}[theorem]{Conjecture}
   
   \newtheorem{problem}[theorem]{Problem}
}
{\theoremstyle{definition}
   
   \newtheorem{example}[theorem]{Example}
   \newtheorem{definition}[theorem]{Definition}
}
{\theoremstyle{remark}
   \newtheorem*{remark}{Remark}
   \newtheorem*{remarks}{Remarks}
}
\newcommand{\RR}{{\mathbb{R}}}

\newcommand{\CC}{{\mathbb{C}}}
\newcommand{\QQ}{{\mathbb{Q}}}

\newcommand{\PP}{{\mathbb{P}}}
\newcommand{\ZZ}{{\mathbb{Z}}}

\newcommand{\bbA}{{\mathbb{A}}}
\newcommand{\bfa}{{\mathbb{A}}}
\newcommand{\bfg}{{\mathbb{G}}}

\newcommand{\cC}{{\mathcal C}}

\newcommand{\cO}{{\mathcal O}}
\newcommand{\cP}{{\mathcal P}}

\newcommand{\Spec}{\operatorname{Spec}}
\newcommand{\Proj}{\operatorname{Proj}}

\newcommand{\Sym}{{\operatorname{Sym}}}

\newcommand{\das}{\dashrightarrow}
\newcommand{\dar}{\downarrow}

\newcommand{\dra}{\dashrightarrow}

\newcommand{\cW}{\cite{Wlodarczyk2}}
\newcommand{\setmin}{\,\protect%
\begin{picture}(8,5)\qbezier(1,4.5)(4,3.)(7,1.5)\end{picture}\,}

\newcommand{\rem}[1]{}
\newcommand{\marg}[1]{}
\begin{document}
\title[Torification and Factorization of Birational Maps]{Torification and 
Factorization of Birational Maps} 

\author{Dan Abramovich}
\thanks{D.A. was partially supported by NSF grant DMS-9700520 and by an
Alfred P. Sloan research fellowship. In addition, he would  like to thank the
Institut des Hautes \'Etudes Scientifiques, Centre Emile Borel (UMS 839,
CNRS/UPMC), and Max Planck Institut f\"ur Mathematik for a fruitful
visiting period.} 
\address{Department of Mathematics\\ Boston University\\ 111 Cummington
Street\\ Boston, MA 02215\\ USA} 
\email{abrmovic@math.bu.edu}
\author{Kalle Karu}
\thanks{K.K. was partially supported by NSF grant DMS-9700520} 
\address{Department of Mathematics\\ Harvard University\\ One Oxford Street \\
Cambridge, MA 02139\\ USA} 
\email{kkaru@math.harvard.edu}
\author{Kenji Matsuki}
\thanks{K.M. has received no financial support from NSF or
NSA during the course of this work.}
\address{Department of Mathematics \\ Purdue University \\  
                        1395 Mathematical Sciences Building \\  West Lafayette,
                        IN 47907-1395 \\ USA}
\email{kmatsuki@math.purdue.edu}
\author{Jaros{\l}aw W{\l}odarczyk}
\thanks{J.W. was supported in part by a Polish KBN grant.}
\address{	Instytut Matematyki UW\\
		Banacha 2, 02-097 Warszawa\\ Poland}
\email{jwlodar@mimuw.edu.pl}
\date{\today}

\begin{abstract}
Building on work of the fourth author in \cite{Wlodarczyk2}, we prove the weak
factorization conjecture for birational maps in characteristic zero: 
a birational map between complete
nonsingular varieties over an algebraically closed field $K$ of characteristic
zero is a composite of blowings up and blowings down with smooth centers.  
\end{abstract}

\maketitle

\tableofcontents
\addtocounter{section}{-1}

\section{Introduction}
We work over an algebraically closed  field  $K$ of characteristic
0. We denote the multiplicative group of $K$ by $K^*$.

\subsection{Statement of the main result}

The purpose of this paper is to give a proof for the following weak
factorization conjecture of birational maps. 
We note that another proof of this theorem
was  given by the fourth author in
\cite{Wlodarczyk3}. 
 See section \ref{Subsec:comparison} for a brief
comparison of the two approaches. 

\begin{theorem}[Weak Factorization]
\label{Th:weak-factorization}
 Let
 $\phi:X_1 \dashrightarrow X_2$ be a birational map between complete
 nonsingular algebraic varieties $X_1$ and $X_2$ over an algebraically closed
 field $K$ of characteristic zero, and let $U\subset  X_1$ be an
 open set where $\phi$ is an isomorphism.  Then $\phi$ can be factored into a
 sequence of blowings up and blowings down with smooth irreducible
 centers disjoint from  $U$, namely, there exists a 
 sequence of birational maps between complete nonsingular algebraic varieties
$$
 X_1 = V_0 \stackrel{\varphi_1}{\dashrightarrow} V_1
 \stackrel{\varphi_2}{\dashrightarrow} \cdots
 \stackrel{\varphi_i}{\dashrightarrow} V_i
 \stackrel{\varphi_{i+1}}{\dashrightarrow} V_{i+1}
 \stackrel{\varphi_{i+2}}{\dashrightarrow}
 \cdots \stackrel{\varphi_{l-1}}{\dashrightarrow}
 V_{l-1} \stackrel{\varphi_l}{\dashrightarrow} V_l = X_2
$$
   where
 \begin{enumerate}
   \item  
   $\phi = \varphi_l \circ \varphi_{l-1} \circ \cdots \varphi_2 \circ
 \varphi_1$,
   \item $\varphi_i$ are isomorphisms on $U$, and 
  \item either 
  $\varphi_i:V_i \dashrightarrow V_{i+1}$ or $\varphi_i^{-1}:V_{i+1}
\dashrightarrow V_{i}$ 
  is a morphism obtained by blowing up a smooth irreducible center
 disjoint from $U$. 
 \end{enumerate}
Furthermore, there is an index $i_0$ such that for all $i\leq i_0$ the map
$V_i\das X_1$ is a projective morphism, and for all $i\geq i_0$ the map
$V_i\das X_2$ is a projective morphism. In particular, if $X_1$ and
$X_2$ are  projective  then all the $V_i$ are  projective. 
\end{theorem}

\subsection{Strong factorization}
If we insist in the assertion above that
$\varphi_1^{-1},\ldots,\varphi_{i_0}^{-1} $ and
$\varphi_{i_0+1},\ldots,\varphi_l$  
be regular maps for some $i_0$, we obtain the following strong
factorization conjecture.    

\begin{conjecture}[Strong Factorization]
Let the situation be as in Theorem \ref{Th:weak-factorization}. Then
there exists a diagram 
$$\begin{array}{rcccl}
 & & Y & & \\
 & \psi_1\swarrow & & \searrow \psi_2 & \\
X_1 & &  \stackrel{\phi}{\dashrightarrow} & & X_2 
\end{array}$$
where the morphisms $\psi_1$ and $\psi_2$ are composites of blowings
up of smooth centers disjoint from $U$.
\end{conjecture}

See Section \ref{Sec:strong-factorization} for further discussion.

\subsection{Generalizations of the main theorem}
We consider the following categories, in which we denote the morphisms by
``broken arrows'':  
\begin{enumerate} 
\item the objects  are  complete nonsingular algebraic
spaces over an arbitrary 
field $L$ of characteristic 0, and  broken
arrows $X \das Y$ denote birational 
$L$-maps, and 
\item the objects  are  compact complex manifolds,  and  broken
arrows $X \das  
Y$ denote bimeromorphic maps. 
\end{enumerate}

Given two broken arrows $\phi: X\das Y$ and
$\phi':X'\das Y'$
 we define an {\em absolute
isomorphism} $g:\phi\to \phi'$
  as follows:
\begin{itemize} 
\item In case $X$ and $Y$ are algebraic spaces over $L$, and $X'$, $Y'$ are
over $L'$, then $g$ consists of an isomorphism $\sigma:\Spec L \to \Spec L'$,
together with a pair of biregular $\sigma$-isomorphisms $g_X : X \to X'$ and
$g_Y: Y \to 
Y'$, such that $\phi'\circ g_X = g_Y\circ
\phi$.
\item In the analytic case, $g$ simply consists of a pair of biregular
isomorphisms $g_X : X \to X'$ and $g_Y: Y \to   
Y'$, such that $\phi'\circ g_X = g_Y\circ
\phi$.
\end{itemize}

\begin{theorem}
\label{Th:general-weak-factorization}
 Let
 $\phi:X_1 \dashrightarrow X_2$ be as in case (1) or (2) above. Let
 $U\subset  X_1$ be an 
 open set where $\phi$ is an isomorphism. Then $\phi$ can be factored, 
 functorially with respect to absolute isomorphisms, into a
 sequence of  blowings up and blowings down with smooth
 centers disjoint from  $U$.  Namely, to any such $\phi$ we associate a
 diagram in  the  corresponding category
$$
 X_1 = V_0 \stackrel{\varphi_1}{\dashrightarrow} V_1
 \stackrel{\varphi_2}{\dashrightarrow} \cdots
 \stackrel{\varphi_i}{\dashrightarrow} V_i
 \stackrel{\varphi_{i+1}}{\dashrightarrow} V_{i+1}
 \stackrel{\varphi_{i+2}}{\dashrightarrow}
 \cdots \stackrel{\varphi_{l-1}}{\dashrightarrow}
 V_{l-1} \stackrel{\varphi_l}{\dashrightarrow} V_l = X_2
$$
   where
 \begin{enumerate}
   \item  
   $\phi = \varphi_l \circ \varphi_{l-1} \circ \cdots \varphi_2 \circ
 \varphi_1$,
   \item $\varphi_i$ are isomorphisms on $U$, and 
  \item either 
  $\varphi_i:V_i \dashrightarrow V_{i+1}$ or $\varphi_i^{-1}:V_{i+1}
\dashrightarrow V_{i}$ 
  is a morphism obtained by blowing up a smooth center
 disjoint from $U$.
 \item Functoriality: if $g:\phi\to \phi'$ is an absolute isomorphism, carrying
 $U$ to $U'$, and
 $\varphi_i':V_i' 
 \dashrightarrow V_{i+1}'$ is the factorization of $\phi'$, then the resulting
 rational maps $g_i:V_i\das V_i'$ are biregular.
 \item Moreover,
   there is an index $i_0$ such that for all $i\leq i_0$ the map 
   $V_i\das X_1$ is a projective morphism, and for all $i\geq i_0$ the map
  $V_i\das X_2$ is a projective morphism.
 \item Let $E_i\subset V_i$ be the exceptional divisor of  $V_i \to X_1$
 (respectively, $V_i \to X_2$) in case $i\leq i_0$ (respectively, $i\geq
 i_0$). Then the above centers of blowings up in $V_i$ have normal crossings
 with  $E_i$. If, moreover, $X_1\setmin U$ (respectively, $X_2\setmin U$) is a
 normal crossings divisor, then the centers of blowing up have normal crossings
 with the inverse images of this divisor.
 \end{enumerate}

\end{theorem}

\begin{remarks}
\begin{enumerate}
\item Note that, in order to achieve functoriality, we cannot require the
centers of blowing up to be irreducible. 
\item Functoriality implies, as immediate corollaries, the existence of
factorization over 
any field of characteristic 0, as well as factorization, equivariant under the
action of a group $G$, of
a $G$-equivariant birational map.
\item The same theorem holds true for varieties or algebraic
spaces of dimension $d$ over a perfect  field of characteristic $p>0$ 
{\em assuming that canonical embedded resolution of singularities holds  true
for varieties or 
algebraic spaces of dimension $d+1$ in characteristic $p$.} The proof for
varieties goes through word for word as in this paper, while for the algebraic
space case one needs to recast some of our steps from the Zariski topology to
the \'etale topology (see \cite{Kato}, \cite{Matsuki-notes}).
\item
While this theorem clearly implies the main theorem as a special case, we
prefer to carry out the proof of the main theorem throughout the text, and 
indicating the changes one needs to perform
for proving Theorem \ref{Th:general-weak-factorization} in section
\ref{Sec:generalizations}.      
\item This is by no means the most general case to which our methods apply, and
we are aware of some applications which are not covered in this
statement. It may be of interest in the future to codify a minimal set of
axioms needed to carry out this line of proof of weak factorization.
\end{enumerate}
\end{remarks}

\subsection{Early origins of the problem}
The history of the factorization problem of birational maps could  be
traced back to the Italian school of algebraic geometers, who already
knew that the operation of blowing up points on 
surfaces is a fundamental  source of richness for surface  geometry: the
importance of the strong factorization theorem in dimension 2 (see
\cite{Zariski}) 
cannot be overestimated in the analysis of the birational geometry of algebraic
surfaces. We can only guess that Zariski, possibly even members of the
Italian school, contemplated the problem in higher dimension early on,
but refrained from
stating it before results on  resolution of singularities were
available.  The question of strong factorization  was explicitly stated
by Hironaka as  ``Question (F$'$)'' in \cite{Hironaka3}, Chapter 0,
\S6, and the question of weak factorization was raised in
\cite{Miyake-Oda}. The problem remained largely open in 
higher dimensions despite the efforts and interesting results of many
(see e.g. Crauder \cite{Crauder}, Kulikov \cite{Kulikov}, Moishezon
\cite{Moishezon}, Schaps
\cite{Schaps}, Teicher \cite{Teicher}). Many of these were summarized  by
Pinkham \cite{Pinkham}, where  the weak factorization conjecture is explicitly
stated. 

\subsection{The toric case}
For toric birational maps,  the equivariant versions of
the weak and strong factorization conjectures were posed in
\cite{Miyake-Oda} and came to be known as Oda's weak and strong conjectures.
While the toric version can be viewed as a special case of the 
general factorization conjectures,  many of the examples
demonstrating the difficulties  in higher dimensions are in fact toric (see
Hironaka \cite{Hironaka2}, Sally \cite{Sally}, Shannon
\cite{Shannon}).  Thus Oda's conjecture presented a substantial challenge and 
combinatorial difficulty. In dimension 3, Danilov's proof of Oda's weak
conjecture \cite{Danilov2} was
later supplemented by Ewald \cite{Ewald}. Oda's weak conjecture was solved in
arbitrary dimension by  J. W{\l}odarczyk in  
\cite{Wlodarczyk1}, and another proof was given by R. Morelli in
\cite{Morelli1} (see also \cite{Morelli2},  
\cite{Abramovich-Matsuki-Rashid}). An important combinatorial notion
which Morelli introduced into this study is that of a {\em cobordism} between
fans. The algebro-geometric realization of Morelli's combinatorial
cobordism is the notion of a {\em birational cobordism} introduced in \cW.

In \cite{Morelli1}, R. Morelli also proposed a proof of Oda's {\em strong}
conjecture. A gap in this proof, which was not noticed in
\cite{Abramovich-Matsuki-Rashid}, was recently discovered by K.\ Karu. As far
as we know, Oda's strong conjecture stands unproven at present even in
dimension 3.

\subsection{The local version}
There is a local version of the factorization conjecture, formulated and
proved in dimension 2 by Abhyankar (\cite{Abhyankar}, Theorem 3).
Christensen \cite{Christensen} posed the problem in general and solved it for
some special cases in dimension 3. Here the 
varieties $X_1$ and $X_2$ are replaced by appropriate birational local rings
dominated by  a fixed valuation, and
blowings up are replaced by monoidal transforms subordinate to the 
valuation. The weak form of this local conjecture was recently solved by
S. D. Cutkosky in a series of papers
\cite{Cutkosky1,Cutkosky3}. Cutkosky also shows that the strong version of the
conjecture follows from Oda's strong factorization conjecture for toric
morphisms. In a sense, Cutkosky's result  
says that the only {\em local} obstructions to solving the global strong
factorization conjecture lie in the toric case.

\subsection{Birational cobordisms}
Our method is based upon the theory of birational
cobordisms \cite{Wlodarczyk2}. As mentioned above, this theory was
inspired by the 
combinatorial notion of 
polyhedral cobordisms of R. Morelli \cite{Morelli1}, which was used in his
proof of weak  factorization for toric birational  maps.

Given a birational map
$\phi:X_1 \dashrightarrow X_2$, a {\em birational
cobordism} $B_{\phi}(X_1,X_2)$ is a variety of  dimension $\dim(X_1)+1$ with an
action of the multiplicative group $K^*$.  It is analogous to the
usual cobordism
$B(M_1,M_2)$ between differentiable manifolds $M_1$ and $M_2$ given by a Morse
function $f$. In the differential setting one can construct an action
of the additive real group $\RR$, where the ``time'' $t\in \RR$ acts
as a diffeomorphism induced by integrating the vector 
field $grad(f)$; hence the multiplicative group $(\RR_{>0},\times)=
\exp(\RR,+)$ 
acts as well.  The critical points of $f$ are
precisely the fixed points of the action of the multiplicative group, and 
the homotopy type of fibers of $f$ changes when we pass through these
critical points. Analogously, in the algebraic setting ``passing through''
the fixed points of the $K^*$-action induces a birational
transformation. Looking at the action on the tangent space at each 
fixed point, we obtain a locally toric description of the
transformation.  This already gives the main result of \cW: a factorization  
of $\phi$ into {\em locally toric birational maps} among varieties
with locally toric structures.  Such birational transformations can also be
interpreted using the work of Brion-Procesi, Thaddeus, Dolgachev-Hu and others
(see \cite{Brion-Procesi,Thaddeus1,Thaddeus2,Dolgachev-Hu}), which
describes the change 
of  Geometric Invariant Theory quotient associated to a
change of linearization.
%\footnote{Indeed, the third author of the
%present paper 
%was working on interpreting of Morelli's combinatorial cobordisms
%in terms of geometric invariant theory at the time J. \Wlo\ sent us his
%preprint}.

\subsection{Locally toric versus toroidal structures}
Considering the fact that weak factorization has been proven for
{\em toroidal} birational maps (\cite{Wlodarczyk1}, \cite{Morelli1},
\cite{Abramovich-Matsuki-Rashid}), 
one might na{\"\i}vely think that  a locally toric factorization, as  indicated
in 
the previous paragraph,
would already provide a proof for Theorem  \ref{Th:weak-factorization}.

However, in the locally toric structure obtained from a cobordism, the
embedded tori chosen may vary 
from point to point, while a toroidal structure (see below) requires
the embedded tori to 
be induced from one fixed open set. Thus there is still a
gap between the notion of locally toric birational maps and that of
toroidal birational maps. 

\subsection{Torification}
In order to bridge over this gap, we follow ideas introduced by Abra\-mo\-vich
and De Jong in
\cite{Abramovich-de-Jong}, and blow up suitable open subsets, called {\em
quasi-elementary cobordisms},  of the birational
cobordism 
$B_{\phi}(X_1,X_2)$ 
along {\em torific ideals}.  This operation  induces a toroidal
structure in a 
neighborhood of each connected component
$F$ of the fixed point set, on 
which the action of
$K^*$ is a {\em toroidal action} (we say that the blowing up {\em
torifies} the action of $K^*$).
 Now the birational transformation ``passing
through $F$'' is toroidal. We use  
canonical resolution of singularities to desingularize the resulting varieties,
 bringing ourselves to a situation
where we can apply the factorization theorem for toroidal birational 
maps. This completes the
proof of Theorem  \ref{Th:weak-factorization}.

\subsection{Relation with the minimal model program}
It is worthwhile to note the relation of the factorization problem to the
development of  Mori's program. Hironaka \cite{Hironaka1} used the
cone of effective 
curves to study the properties of birational morphisms.  This direction was
further developed and given a decisive impact by Mori \cite{Mori1},
who
%, motivated,  among other
%things, by Kleiman's criterion for ampleness \cite{Kleiman}, 
introduced the notion of extremal rays and 
systematically used it in an attempt to construct  minimal models in
higher dimension, called {\em the minimal model program}. Danilov
\cite{Danilov2} introduced the notion of {\em canonical and terminal
singularities} in conjunction with the toric factorization problem. This was
 developed by Reid into a  general theory of these singularities
\cite{Reid1,Reid2}, which  appear in an essential way in the minimal model
program. The minimal model program is  so far proven up  
to dimension 3 (\cite{Mori2}, see also
\cite{Kawamata1,Kawamata2,Kawamata3,Kollar1,Shokurov}), and for toric 
varieties in arbitrary dimension (See \cite{Reid3}). In the steps of the
minimal model program one is only allowed to  contract a divisor into a
variety with 
terminal singularities, or to perform a flip, modifying some  codimension $\geq
2$ 
loci. This allows a factorization of a given birational morphism into
such ``elementary operations''. An algorithm to factor birational
{\em maps} among  uniruled varieties, 
known as {\em Sarkisov's program,} has been developed  and carried out in
dimension 3  (see \cite{Sarkisov,Reid4,Corti}, and see  \cite{Matsuki}
for the toric 
case). Still,  we do not know of a way to solve
the classical factorization problem using such a factorization.
 
\subsection{Relation with the toroidalization problem}
In \cite{Abramovich-Karu}, Theorem 2.1, it is proven
 that given 
a  morphism of projective varieties $X \to B$, there are modifications
$m_X:X' \to X$ and  $m_B:B' \to B$, with a lifting $X'\to B'$ which has
a toroidal structure. The {\em toroidalization problem} (see
\cite{Abramovich-Karu}, \cite{Abramovich-Matsuki-Rashid}, \cite{King1})
is that of 
obtaining such $m_X$ and $m_B$ which are composites of blowings up with
smooth centers (maybe even with centers supported only over the locus where
$X \to B$ is not toroidal). 

The proof in \cite{Abramovich-Karu} relies on the work of De Jong
\cite{de Jong} and methods of \cite{Abramovich-de-Jong}. The authors
of the present paper have tried to use these methods to approach the
factorization conjectures, so far without success; one notion we do use
in this paper  is  the torific ideal of
\cite{Abramovich-de-Jong}. It would be interesting if one could turn
this approach on its 
head and prove a result on toroidalization using  factorization.

\subsection{Relation with the proof in
\cite{Wlodarczyk3}}\label{Subsec:comparison} Another proof of  
the weak factorization theorem was given independently by the fourth author
in \cite{Wlodarczyk3}. The main difference of the two approaches is that in
the current paper we are using objects such as torific ideals
defined {\em locally} on each quasi-elementary piece of a cobordism. The
blowing up 
of 
a torific ideal gives the quasi-elementary cobordism a toroidal
structure. These toroidal modifications are then pieced together using
canonical  resolution of singularities. In \cite{Wlodarczyk3} one works {\em
globally}: a combinatorial description of {\em stratified toroidal varieties}
and 
appropriate morphisms between them is given, which allows one to apply
Morelli's $\pi$-desingularization algorithm directly to the entire 
birational cobordism.  The  structure of stratified toroidal variety
 on the cobordism is somewhere in between our notions
of locally toric and toroidal structures.

\subsection{Outline of the paper} 
In section \ref{Sec:prelims} we discuss locally toric
and toroidal structures, and  reduce the proof of Theorem
\ref{Th:weak-factorization} to the case where $\phi$ is a projective birational
morphism.

Suppose now we have a projective birational
morphism $\phi:X_1 \to X_2$.
In section \ref{Sec:cobordism} we apply the theory of
birational cobordisms to obtain a factorization into locally toric
birational maps.  Our cobordism $B$ is 
relatively projective over $X_2$, and using a geometric invariant 
theory analysis, inspired by Thaddeus's work,  
 we show that the intermediate varieties
can be chosen to be projective over $X_2$.

In section \ref{Sec:torification} we utilize a factorization of the cobordism
$B$ into quasi-elementary pieces $B_{a_i}$, and for each piece 
 construct an ideal sheaf $I$
whose blowing up {\em torifies} 
the action of $K^*$ on $B_{a_i}$. In other words, $K^*$ acts toroidally
on the  variety obtained by blowing up $B_{a_i}$ along  $I$.

In section \ref{Sec:connecting} we prove the weak factorization theorem by
putting together the toroidal birational transforms induced by the
quasi-elementary cobordisms.
This is done using
canonical resolution of singularities. 

In section \ref{Sec:generalizations} we prove Theorem
\ref{Th:general-weak-factorization}. We then discuss some problems
related to strong factorization in section \ref{Sec:problems}.

\section{Preliminaries}\label{Sec:prelims}

\subsection{Quotients}
Suppose a reductive group $G$ acts on an algebraic variety $X$. We denote by
$X/G$  the space of orbits,
and by  $X /\!/ G$  the space of equivalence classes of orbits, where the
equivalence relation is generated by the condition that two 
orbits are equivalent if  their closures intersect; such a space is endowed
with 
a  scheme structure which satisfies the usual universal property, if such a
structure 
exists. 

In this paper we will only consider $X /\!/ G$ in situations where the
closure of any orbit contains a unique closed orbit (see Definition
\ref{Def:quasi-elementary}).  Moreover, the quotient morphism $X \to X /\!/ G$
will be affine. When this holds we say that the action of $G$ on $X$ is {\em
relatively affine}.

\subsection{Canonical resolution of singularities and canonical
principalization} \label{Sec:resolution}
In the following, we will use canonical versions of Hironaka's theorems on
resolution of 
singularities and principalization of an ideal, proved in
\cite{Bierstone-Milman,Villamayor}.

\subsubsection{Canonical resolution}
A canonical embedded  resolution of singularities $\widetilde{W}\to W$
is a desingularization procedure consisting
of a composite of blowings up with smooth centers, satisfying a number of
conditions.  In particular 
\begin{enumerate} 
\item ``embedded'' means the following: assume the sequence of blowings up is
applied when 
$W\subset U$ is a closed embedding with
 $U$ nonsingular. Denote by $E_i$
the exceptional divisor at some stage of the blowing up. Then  {\bf (a)} $E_i$
is a normal crossings divisor, and 
has normal crossings with the center of blowing up, and {\bf (b)} at the last
stage $\widetilde{W}$  has normal 
crossings with  $E_i$.
\item ``Canonical'' means ``functorial with respect to smooth maps'', namely,
if $\theta:V\to W$ is a smooth morphism then the 
ideals 
blown up are invariant under pulling back by 
$\theta$; hence $\theta$ can be lifted to a smooth morphism
$\widetilde{\theta}:\widetilde{V}\to 
\widetilde{W}$.
\end{enumerate}
 In particular: {\bf (a)} if $\theta:W \to W$ is an
automorphism (of schemes, not necessarily over $K$) then it can be lifted to
an automorphism $\widetilde{W}\to 
\widetilde{W}$, and {\bf (b)}  the canonical resolution behaves well with
respect to \'etale morphisms: if $V\to W$ is \'etale, we get an \'etale
morphism of canonical resolutions $\widetilde{V}\to\widetilde{W}$. 

An important consequence of these conditions is that {\em all the centers of
blowing up lie over the singular locus of $W$.}

\subsubsection{Compatibility with a normal crossings divisor}
If $W \subset U$ is embedded in a nonsingular variety, and $D\subset U$ is a
normal crossings divisor, then a variant of the resolution procedure allows one
to choose the centers of blowing up to have normal crossings with $D_i+E_i$,
where $D_i$ is the inverse image of $D$. This follows since the resolution
setup, as in \cite{Bierstone-Milman}, allows including such a divisor in ``year
0''. 

\subsubsection{Principalization}
By {\em canonical principalization  of an ideal sheaf} in a 
nonsingular variety we mean ``the canonical embedded resolution of
singularities of the subscheme defined by the ideal sheaf making it a
divisor with normal crossings''; i.e., a composite
of blowings up with smooth centers such that the total
transform of the ideal is a divisor with simple normal crossings. Canonical
embedded 
resolution of singularities of an arbitrary subscheme, not necessarily
reduced or irreducible,  is discussed in  
Section~11 of \cite{Bierstone-Milman}, and this implies canonical
principalization, as 
one simply needs to blow up $\widetilde{W}$ at the last step.

\subsubsection{Elimination of indeterminacies}
Now let $\phi:W_1\dra W_2$ be a proper birational map between nonsingular
varieties, and $U\subset W_1$ an open set on which $\phi$ restricts to an
isomorphism. By {\em elimination of indeterminacies} of $\phi$ we mean a
morphism $e:W_1'\to W_1$,  
obtained by a sequence of blowings up with smooth centers disjoint from $U$,
such that the birational map $\phi\circ e$ is a morphism. 

Elimination of
indeterminacies can be reduced to principalization of an ideal sheaf as
follows.  

We may assume that $\phi^{-1}$ is a morphism; otherwise we replace $W_2$ by
the closure of the graph of $\phi$. 
%If $\phi^{-1}$ is a projective morphism,
%i.e., a blowing up of an ideal sheaf $I$, then we can take for $W_1'$ the
%canonical principalization of $I$. 
%If $\phi^{-1}$ is not projective, 
Now we use
Chow's lemma (Corollary~2, p. 504, \cite{Hironaka4}): there exists an ideal
sheaf $I$ on $W_1$ such that the blowing up of $W_1$ along $I$ factors through
$W_2$. Hence the canonical principalization of $I$ also factors through $W_2$.

Although it is not explicitly stated by Hironaka, the ideal $I$ is
supported in the complement of the open set $U$: the blowing up of $I$ consists
of a sequence of {\em permissible} blowings up (Definition 4.4.3, p. 537,
\cite{Hironaka4}), each of which is supported in the complement of $U$.
Another important fact is, that the ideal $I$ is {\em invariant}, namely, it is
functorial under absolute isomorphisms: if $\phi:W_1'\dra W_2'$ is another
proper birational map, with corresponding ideal $I'$, and  
$\theta_i:W_i \to W_i'$ are  isomorphisms such that  $\phi'\circ \theta_1 =
\theta_2\circ \phi$, then $\theta_1^*I' = I$. This follows
simply because at no point in Hironaka's flattening procedure there is a need
for any choice. 

The same results hold for analytic and algebraic spaces. While Hironaka states
his result only in the analytic setting, the arguments hold in the algebraic
setting as well. See \cite{Raynaud-Gruson} for an earlier 
treatment of the case of varieties.

\subsection{Reduction to projective morphisms}
We start with a birational map 
$$\phi:X_1 \dashrightarrow X_2$$
between complete nonsingular   algebraic varieties $X_1$ and $X_2$
defined over $K$ and restricting to an isomorphism on an open set $U$.   

\begin{lemma}[Hironaka] \label{lem:red-proj}
There is a commutative diagram 
$$\begin{array}{rcl} 
X_1' & \stackrel{\phi'}{\to} & X_2' \\
g_1 \dar & & \dar g_2 \\
X_1 & \stackrel{\phi}{\das} & X_2 
\end{array}$$
such that $g_1$ and $g_2$ are composites of blowings up with smooth
centers disjoint from $U$, and $\phi'$ is a projective birational morphism.
\end{lemma}

{\bf Proof.} By Hironaka's theorem on elimination of indeterminacies,
there is a morphism $g_2:X_2' \to X_2$ which is a composite of blowings up
with smooth centers disjoint from $U$, such that the birational map
$h:=\phi^{-1}\circ g_2 : X_2' \to X_1$ is a morphism:
$$\begin{array}{rcl} 
 &  & X_2' \\
 & h \swarrow & \dar g_2 \\
X_1 & \stackrel{\phi}{\das} & X_2 
\end{array}.$$ 
By the same
theorem, there is a morphism $g_1:X_1' \to X_1$ which is a composite of
blowings up 
with smooth centers disjoint from $U$, such that $\phi':=h^{-1}\circ
g_1: X_1' \to X_2'$ is a morphism. Since the composite $h\circ \phi' = g_1$
is projective, it follows that $\phi'$ is
projective. \qed

Thus we may replace $X_1 \das X_2$ by $X_1' \to X_2'$ and assume from
now on that $\phi$ is a projective morphism.

Note that, by the properties of principalization and Hironaka's flattening, the
formation of $\phi':X_1' \to X_2'$ is functorial under absolute isomorphisms,
and the blowings up have normal crossings with the approapriate divisors. This
will be used in the proof of Theorem \ref{Th:general-weak-factorization}.

\subsection{Toric varieties} Let $N\cong\ZZ^n$ be a lattice and
$\sigma\subset N_\RR$ a strictly convex rational polyhedral cone. We denote
the dual lattice by $M$ and the dual cone by $\sigma^\vee\subset M_\RR$. The
{\em affine toric variety} $X=X(N,\sigma)$ is defined as
\[ X=\Spec K[M\cap\sigma^\vee].\]
For  $m\in M\cap\sigma^\vee$ we denote its image in the semigroup algebra
$K[M\cap\sigma^\vee]$ by $z^m$. 

More generally, the toric variety corresponding to a fan $\Sigma$ in $N_\RR$ is
denoted by $X(N,\Sigma)$.

If $X_1=X(N,\Sigma_1)$ and $X_2=X(N,\Sigma_2)$ are two toric varieties, the
embeddings of the torus $T=\Spec K[M]$ in both of them define a 
toric (i.e., $T$-equivariant) birational  map $X_1\dra X_2$.

Suppose  $K^*$ acts effectively on an affine toric variety $X=X(N,\sigma)$ as a
one-parameter subgroup of the torus $T$, corresponding to a primitive lattice 
point $a\in N$. If $t\in K^*$ and $m\in M$, the action on the monomial $z^m$
is given by 
\[ t^*(z^m) = t^{(a,m)}\cdot z^m,\]
where $(\cdot,\cdot)$ is the natural pairing on $N\times M$. The
$K^*$-invariant monomials correspond to the lattice points $M\cap a^\perp$,
hence 
\[ X/\!/K^* \cong \Spec K[M\cap\sigma^\vee\cap a^\perp].\]
If $a\notin \pm\sigma$ then $\sigma^\vee\cap a^\perp$ is a
full-dimensional cone in $a^\perp$, and it follows that $X/\!/K^*$
is again an affine toric variety, defined by the lattice $\pi(N)$ and
cone $\pi(\sigma)$, where $\pi: N_\RR\to N_\RR/\RR\cdot a$ is the
projection. This quotient is a geometric quotient precisely when $\pi:
\sigma\to\pi(\sigma)$ is a bijection.

\subsection{Locally toric and toroidal structures} There is some confusion
in the literature between the notion of {\em toroidal embeddings} and
toroidal morphisms (\cite{KKMS}, \cite{Abramovich-Karu}) and that of {\em
toroidal varieties} (see \cite{Danilov1}), which we prefer to call
{\em locally toric varieties}, and locally toric morphisms between them. A
crucial issue in this paper is the distinction between the two notions.

\begin{definition} 
\begin{enumerate}
\item A variety $W$ is {\em locally toric} if for every closed point $p\in W$
there exists an open neighborhood $V_p\subset W$ of $p$ and an \'etale
morphism $\eta_p:V_p\to X_p$ to a toric variety $X_p$. Such a morphism
$\eta_p$ is called a {\em toric chart} at $p$.
\item An open embedding $U\subset W$ is a {\em toroidal embedding} if  for
every closed point $p\in W$ there exists a toric chart $\eta_p:V_p\to X_p$ at
$p$ such that $U\cap V_p = \eta_p^{-1}(T)$, where $T\subset X_p$ is the 
torus. We call such charts {\em toroidal}. Sometimes we omit the open set $U$
from the notation and simply say that a variety is toroidal.
\item We say that a locally toric (respectively, toroidal) chart on a
variety is compatible with a divisor $D\subset W$ if  $\eta_p^{-1}(T)
\cap D = \emptyset$, i.e., $D$ corresponds to a toric divisor on $X_p$.
\end{enumerate} 
\end{definition}

\begin{definition}\label{Def:proper-birational-morphism}
\begin{enumerate}
\item A proper birational morphism of locally toric varieties $f: W_1 \to
W_2$ is said 
to be {\em locally toric} if for every closed point $q\in W_2$, and any $p\in
f^{-1}q$,  there is a
diagram of {\em fiber squares}
$$\begin{array}{rcccl}
X_p & \leftarrow & V_p & \subset & W_1 \\
\phi\dar & & \dar & & \dar f \\
X_{q} & \leftarrow & V_{q} & \subset & W_2
\end{array}$$
where 
\begin{itemize} 
\item $\eta_p:V_p\to X_p$ is a toric chart at $p$,
\item $\eta_{q}:V_{q}\to X_{q}$ is a toric chart at $q$, and
\item $\phi: X_p\to X_{q}$ is a toric morphism.
\end{itemize}

\item Let $U_i \subset W_i$ $(i=1,2)$ be toroidal embeddings. A   proper
birational 
morphism $f: W_1 \to W_2$ is said to be {\em toroidal}, if it satisfies the
condition above for being locally toric, with toroidal charts. In particular,
$f^{-1}(U_2) = U_1$.
\end{enumerate}
\end{definition}

\begin{remarks}
\begin{enumerate}
\item
A toroidal embedding as defined above is a toroidal embedding without
self-intersection according to the definition in \cite{KKMS}, and a
birational toroidal morphism satisfies the condition of allowability in
\cite{KKMS}. 
\item We note that this definition of a toroidal morphism, where the charts are
chosen locally on the target $W_2$, {\em differs}
from that in \cite{Abramovich-Karu}, where the charts are taken locally in the
source $W_1$. It is a nontrivial fact, which we will not need in this
paper, that {\em these notions do agree for proper birational morphisms.} 
\item As the reader may  notice, 
one can propose several variants of the definition above as well as the ones to
come, and one can raise many subtle questions about comparison between
the resulting notions. We will not  address these issues in this paper at all. 
\item
To a toroidal embedding $(U_W\subset W)$ one can associate a
polyhedral complex $\Delta_W$, such that proper birational toroidal morphisms
to $W$, 
up to isomorphisms, 
are in one-to-one correspondence with certain subdivisions of the complex
(see \cite{KKMS}). It follows from this that the composition of two proper
birational toroidal morphisms $W_1\to W_2$ and $W_2\to W_3$ is again
toroidal: the first morphism corresponds to a subdivision of $\Delta_{W_2}$,
the second one to a subdivision of $\Delta_{W_3}$, hence their composition is
the unique toroidal morphism correponding to the subdivision $\Delta_{W_1}$
of $\Delta_{W_3}$.
\item
A composition of locally toric birational morphisms is not locally toric in
general. A simple example is given by blowing up a point on a
nonsingular threefold, and then blowing up a nonsingular curve tangent to the
exceptional divisor. 
\item
One can define notions of locally toric and toroidal morphisms even when
the morphisms are not proper or birational (see,
e.g., \cite{Kato}, \cite{Abramovich-Karu}). We 
will not need such notions in this paper.
\item
Some of the issues we avoided discussing here are addressed in the third
author's lecture notes \cite{Matsuki-notes}.
\end{enumerate}
\end{remarks}

\begin{definition}[\cite{Hironaka3},\cite{Iitaka}] Let $\psi : W_1
\das W_2$ be a rational map defined on a dense open
subset $U$. Denote by $\Gamma_\psi$ the closure of the graph of  $\psi_U$
in $W_1 \times W_2$. We say that $\psi$ is
{\em proper} if the projections $\Gamma_\psi\to W_1$ and
$\Gamma_\psi\to W_2$ are both proper.
\end{definition}

\begin{definition} \label{def:map-loc-tor}
\begin{enumerate}
\item A proper birational map $\psi:W_1 \dashrightarrow W_2$
between two locally toric varieties $W_1$ and $W_2$ is said to be {\em
locally toric} if there exists a locally toric variety $Z$ and a commutative
diagram 
$$\begin{array}{ccccc}
 & & Z & & \\
 & \swarrow & & \searrow & \\
W_1 & &\stackrel{\psi}{\das} & & W_2
\end{array}
$$
where $Z \to W_i$ $(i=1,2)$ are proper birational locally toric morphisms.
\item Let $U_i \subset W_i$ be toroidal embeddings. A proper birational map
$\psi:W_1 \dashrightarrow W_2$ is said to be {\em toroidal} if there exists a
toroidal embedding $U_Z \subset Z$ and a diagram as above where $Z \to W_i$
are proper birational toroidal  morphisms. In particular, a proper birational
toroidal map induces an isomorphism between the open sets $U_1$ and $U_2$.
\end{enumerate}
\end{definition}
\begin{remarks}
\begin{enumerate}
\item
It follows from the correspondence between toroidal modifications and
subdivisions of polyhedral complexes that the composition of toroidal
birational maps given by  $W_1 \leftarrow Z_1 \to W_2$ and  $W_2 \leftarrow Z_2
\to W_3$ is again toroidal. Indeed, if $Z_1\to W_2$ and $Z_2\to W_2$
correspond to two subdivisions of 
$\Delta_{W_2}$, then a common refinement of the two subdivisions corresponds to
a 
toroidal embedding $Z$ such that $Z\to Z_1$ and $Z\to Z_2$ are toroidal
morphisms. For example, the coarsest refinement correpsonds to taking for $Z$
the normalization of the closure of the graph of the birational map 
$Z_1\das Z_2$. The composite maps $Z \to W_i$ are all toroidal birational
morphisms. 
% Also, if a
% morphism between toroidal embeddings is toroidal as a birational map, it is
% toroidal as a morphism. This follows from \cite{Abramovich-Karu}.
\item It can be shown that a toroidal birational map which is regular is a
toroidal morphism, therefore definitions \ref{Def:proper-birational-morphism}
and \ref{def:map-loc-tor} are compatible in the toroidal situation. We do not
know if this is true for locally toric maps.  
\end{enumerate}
\end{remarks}
  
A composition of locally toric birational maps is not locally toric in
general. Even worse, the locally toric
structures on the morphisms $Z\to W_i$ may be given with respect to
incompatible toric charts in $Z$. Hence, for points $p\in W_1$ and $q\in W_2$
there may exist no toric charts at $p$ and $q$ in which the map $\psi$ is
given by a birational toric map. 

To remedy this, we define  a stronger
version of locally toric and toroidal maps. These are the only maps we will
need in the considerations of  the current paper.

\begin{definition}
\begin{enumerate}
\item A proper birational map $\psi:W_1 \dashrightarrow W_2$ between locally
toric varieties $W_1$ and $W_2$ is called {\em tightly locally toric} if
there exists a locally toric variety $Y$ and a commutative a diagram 
$$\begin{array}{ccccc}
W_1 & &\stackrel{\psi}{\das} & & W_2\\
 & \searrow & & \swarrow & \\
 & & Y &, & \\
\end{array}$$
where $W_i \to Y$ are proper birational locally toric morphisms, and for
every closed point
$q\in Y$ there exists a toric chart $\eta_q:V_q\to X_q$ at $q$ such that the
morphisms $W_i\to Y$ can be given locally toric structures with respect to the
same chart $\eta_q$.
\item Let $U_i \subset W_i$ be toroidal embeddings. A proper birational map
$\psi:W_1 \dashrightarrow W_2$ is said to be {\em tightly toroidal} if there
exists a toroidal embedding $U_Y \subset Y$ and a diagram as above where
$W_i\to Y$ are proper birational toroidal  morphisms.
\end{enumerate}
\end{definition}

\begin{remark}
The argument used before to show that a composition of toroidal birational maps
is toroidal, shows that a tightly toroidal map is toroidal. A 
composition of tightly toroidal maps is not tightly toroidal in general. As
for  tightly locally toric maps, all varieties and morphisms can be given
toroidal structures locally in $Y$. Now letting $Z$ be the normalization of
the closure of the graph of $\psi$, it follows that $\psi$ is locally toric. 
\end{remark}

\subsection{Weak factorization for toroidal birational maps}
 The weak factorization theorem for proper
birational toric maps can 
be extended to the case of proper birational toroidal maps. This is
proved in \cite{Abramovich-Matsuki-Rashid} for toroidal morphisms, using the
correspondence between birational toroidal morphisms and subdivisions of
polyhedral 
complexes. 
The general case of a toroidal birational map $W_1\leftarrow Z \to W_2$ can be
deduced 
from this, as follows. By toroidal resolution of singularities we may assume
$Z$ is nonsingular. We apply toroidal weak factorization to the morphisms
$Z\to  W_i$, to get a sequence of toroidal birational maps
$$W_1 = V_1  \das V_2 \das\cdots\das V_{l-1} \das V_l= Z \das V_{l+1}
\das\cdots\das V_{k-1}\das V_k=W_2$$ 
consisting of smooth toroidal blowings up and down.

We state this result for later reference:

\begin{theorem}\label{thm:morelli} Let $U_1\subset W_1$ and
$U_2\subset W_2$ be nonsingular toroidal 
embeddings.
Let $\psi: W_1\dra W_2$ be a proper toroidal birational
map. Then $\phi$ can be factored into a
sequence of toroidal birational maps consisting of smooth toroidal  blowings
up and down.  
\end{theorem}

This does not immediately imply that one can choose a
factorization satisfying a projectivity statement as in the main theorem, or
in a  functorial manner. We  will show these facts in Sections
\ref{Sec:proj-tor-fac} and
\ref{Sec:generalizations}, respectively. It should be mentioned that if toric
{\em strong} 
factorization is true, then the toroidal case follows. 

\subsection{Locally toric and toroidal actions} \label{Sec:loc-toric-action}

\begin{definition}[see \cite{Mumford-Fogarty-Kirwan}, p. 198] Let $V$ and $X$
be affine varieties with $K^*$-actions, and let $\eta: V\to X$ be a
$K^*$-equivariant \'etale morphism. Then 
$\eta$ is said to be {\em strongly \'etale} if
\begin{enumerate}
\item[(i)] the quotient map $V/\!/K^*\to X/\!/K^*$ is \'etale, and 
\item[(ii)] the natural map
\[ V\to X\mathbin{\mathop{\times}\limits_{X/\!/K^*}} V/\!/K^*\]
is an isomorphism.
\end{enumerate}
\end{definition}

\begin{definition} 
\begin{enumerate}
\item Let $W$ be a locally toric variety with a  $K^*$-action, such that $W
/\!/ K^*$ exists.
 We
say that the action is {\em locally toric} if for any closed point $p\in W$ we
have 
a toric chart $\eta_p:V_p\to X_p$ at $p$, and a one-parameter subgroup
$K^*\subset T_p$ of the torus in $X_p$,  satisfying
\begin{itemize}
\item $V_p= \pi^{-1} \pi V_p$, where $\pi:W\to W/\!/K^*$ is the projection;
\item $\eta_p$ is $K^*$-equivariant and strongly \'etale.
\end{itemize}
\item If $U\subset W$ is a toroidal embedding, we say that $K^*$ acts
toroidally on $W$ if the charts above can be chosen toroidal.
\end{enumerate}
\end{definition}

The definition above is equivalent to the existence of the following diagram
of fiber squares:
$$\begin{array}{ccccl}
X_p & \leftarrow & V_p & \subset & W \\
\dar & & \dar & & \dar f \\
X/\!/K^* & \leftarrow & V_p/\!/K^* & \subset & W/\!/K^*
\end{array}$$
where the horizontal maps provide toric (resp. toroidal) charts in $W$ and
$W/\!/K^*$. It follows that the quotient of a locally toric variety by a
locally toric action is again locally toric; the same holds in the toroidal
case. 

\begin{remark}
If we
do not insist on the charts being strongly \'etale, then the morphism of
quotients may fail to be \'etale. Consider, for
instance,  the space $X = \Spec\ K[x, x^{-1},y]$ with the action $t(x,y) =
(t^2x, t^{-1}y)$. The quotient is $X / K^*=\Spec\ K[xy^2].$ There is an
equivariant 
\'etale cover $V =  \Spec\ K[u, u^{-1},y]$ with the action $t(u,y) =
(tu, t^{-1}y)$, where the map is defined by $x = u^2$. The quotient is $V /
K^*=\Spec\ K[uy],$ which is a {\em branched} cover of $X/K^*$, since $xy^2 =
(uy)^2$. 
\end{remark}

The following lemma shows that locally toric $K^*$-actions are ubiquitous. We
note that it can be proven with fewer assumptions, see \cite{Wlodarczyk2},
\cite{Matsuki-notes}. 

\begin{lemma}\label{lem:loc-tor-act} Let $W$ be a nonsingular variety with a
relatively affine $K^*$-action, that is, the scheme $W/\!/K^*$ exists and the
morphism 
$W\to W/\!/K^*$ is an affine morphism. Then the action of $K^*$ on  $W$ is
locally toric.
\end{lemma}

{\bf Proof.}  Taking an affine open in $W/\!/K^*$, we may assume that $W$ is
affine. We embed $W$ equivariantly into a projective space and take its
completion (see, e.g., \cite{Sumihiro}). After applying equivariant resolution
of singularities to this 
completion (see Section \ref{Sec:resolution}) we may also assume that
$\overline{W}$ is a 
nonsingular projective variety with a $K^*$-action, and $W\subset \overline{W}$
is an affine invariant open subset.

Let $p\in W$ be a closed point. Since $\overline{W}$ is complete, the orbit
of $p$ has a limit point $q=\lim_{t\to 0} t(p)$ in $\overline{W}$. Now $q$ is
fixed by $K^*$, hence $K^*$ acts on the cotangent space $m_q/m_q^2$ at
$q$. Since $K^*$ is reductive, we can lift a set of eigenvectors of this action
to semi-invariant local 
parameters $x_1,\ldots, x_n$ at $q$. These local parameters define a
$K^*$-equivariant  \'etale morphism $\eta_q: V_q\to X_q$ from an affine $K^*$
invariant open neighborhood $V_q$ of $q$ to the tangent space $X_q =
\Spec (\Sym\ m_q/m_q^2)$ at $q$. The latter has a structure of a toric variety,
where the torus is the complement of the zero set of $\prod x_i$.  

Separating the parameters $x_i$ into $K^*$-invariants and non-invariants, we
get a factorization $X_q = X_q^0 \times X_q^1$, where the action of $K^*$ on
$X_q^1$ is trivial, and the action on $X_q^0$ has $0$ as its unique fixed
point. Thus we get a product decomposition $X_q/\!/K^* =  X_q^0/\!/K^* \times
X_q^1$. 

By  Luna's Fundamental Lemma (\cite{Luna}, Lemme 3), there exist affine
$K^*$-invariant neighborhoods $V_q'$ of $q$ and $X_q'$ of $0$,  such that
the restriction $\eta_q': V_q' \to X_q'$ is  
strongly \'etale. Consider first the case $q\in W$, in which case we may
replace $p$ by $q$. Denote $Z=X_q^{K^*}\cap
X_q'$. Then $Z \subset X_q^{K^*}\simeq  X_q^1$ is affine open, and,  using the
direct product decomposition above, $X_q^0\times 
Z \subset X_q$ is affine open. 
Denote $X_q'' = X_q'\cap X_q^0\times Z$. This is affine open in $X_q$, and it
is easy to see that $X_q''/\!/K^* \to X_q/\!/K^*$ is an open embedding: an
orbit in $X_q''$ is closed if and only if it is closed in $X_q$. Writing 
$V_q'' = {\eta_q'}^{-1}X_q''$, it follows that $V_q'' \to X_q$ is a strongly
\'etale toric chart.

In case $q\notin W$, replace $V_q$ by $V_q''$ . Now $\eta_q$
is injective on any orbit, and 
therefore it is injective on the orbit of $p$. Let $X_p\subset X_q$ be the
affine open toric subvariety in which the torus orbit of $\eta_q(p)$ is closed,
and 
let $V_p = \eta_q^{-1}X_p\cap W$.  Now consider the restriction $\eta: V_p\to
X_p$, where the $K^*$-orbits of $p$ and $\eta(p)$ are closed. By Luna's
Fundamental Lemma  there exist 
affine open 
$K^*$-invariant neighborhoods $V_p'\subset V_p$  and $X_p'\subset X_p$ of
$\eta_p(p)$ 
such that the restriction $\eta: V_p'\to X_p'$ is a strongly \'etale
morphism. Since $X_p/K^*$ is a geometric quotient, we may assume $X_p'=X_p$ and
we have a strongly \'etale toric chart.

It remains to show that the charts can be chosen saturated with respect to
the projection $\pi: W\to W/\!/K^*$. If the orbit of $p$ has a limit point
$q=\lim_{t\to 0} t\cdot p$ or $q=\lim_{t\to \infty} 
t\cdot p$ in $W$, which is necessarily unique as $\pi$ is affine, then an
equivariant  toric
chart at $q$  
also covers $p$. So we may replace $p$ by $q$ and assume that the orbit of $p$
is closed. Now $\pi (W\setmin V_p)$ is closed and does
not contain $\pi(p)$, so we can choose an affine neighborhood $Y$ in its
complement, and replace $V_p$ by $\pi^{-1} Y$.\qed

\section{Birational Cobordisms}\label{Sec:cobordism}
\subsection{Definitions}

\begin{definition}[\cW] \label{Def:cobordism}
Let $\phi:X_1 \dashrightarrow X_2$ be a birational
map between two algebraic varieties $X_1$ and $X_2$ over $K$, isomorphic on an
open set $U$.  A normal algebraic variety
$B$ is called a {\em birational cobordism} for $\phi$ and denoted by
$B_{\phi}(X_1,X_2)$ if it satisfies the following conditions.
\begin{enumerate}
\item The multiplicative group $K^*$ acts effectively on $B =
B_{\phi}(X_1,X_2)$. 
\item The sets
$$\begin{array}{rcrl}
B_- &:=& \{x \in B:&\lim_{t \to 0\ }t(x) \mbox{ does not exist in }B\} \\ 
\mbox{and  } 
B_+ &:=& \{x \in B:&\lim_{t \to \infty}t(x) \mbox{ does not exist
in }B\} 
\end{array}
$$
are nonempty Zariski open subsets of $B$. 
\item There are isomorphisms
$$ B_-/K^* \stackrel{\sim}{\to} X_1 \quad \mbox{ and } \quad
B_+/K^* \stackrel{\sim}{\to} X_2.  
$$
\item Considering the rational map $\psi:B_-\das B_+$ induced
by the inclusions $B_-\cap B_+\subset B_-$ and $B_-\cap B_+\subset
B_+$, the following diagram commutes:
$$
\begin{array}{ccc} B_- & \stackrel{\psi}{\das}& B_+ \\
  \dar & & \dar \\ X_1 &\stackrel{\phi}{\das} & X_2\end{array}.
$$
\end{enumerate}
We say that $B$ {\em respects the open set $U$} if $U$ is contained in
the image of $(B_-\cap B_+)/K^*$.
\end{definition} 

%\begin{remark} The definition implies that the action of $K^*$ on $B$ is
%nontrivial. If needed, we can always reduce to the case where the action is
%effective.  
%\end{remark}

\begin{definition}[\cW]  Let $B = B_\phi(X_1,X_2)$ be a
birational cobordism, 
and let $F\subset B^{K^*}$ be a subset of the fixed-point
set.
We define 
$$\begin{array}{rcl} 
F^+ & = & \{ x\in B | \lim_{t\to 0\ } t(x)\in F\} \\ 
 F^- & = & \{ x\in B | \lim_{t\to \infty}t(x)\in F\} \\ 
F^\pm & = & F^+ \cup  F^- \\
F^*\ & = & F^{\pm} \setmin F
\end{array}$$
\end{definition}

\begin{definition}[\cW]
Let $B=B_\phi(X_1,X_2)$ be a  birational cobordism. We define a
relation $\prec$ among  {\em connected components} of $B^{K^*}$ as follows:
let $F_1, F_2\subset B^{K^*}$ be  two connected components, and set $F_1
\prec F_2$ if there is a point $x\notin B^{K^*}$ such that $\lim_{t\to 0 } t(x)
\in F_1$ and  $\lim_{t\to \infty}t(x)\in F_2$. 
\end{definition}

\begin{definition} \label{Def:quasi-elementary} A  birational cobordism
$B=B_\phi(X_1,X_2)$ is said to be {\em quasi-elementary} if any two
connected components $F_1, F_2\subset B^{K^*}$ are incomparable with respect
to $\prec$.
\end{definition}

Note that this condition prohibits, in particular, the existence of a ``loop'',
namely a connected 
component $F$ and a point $y\notin F$ such that both $\lim_{t\to 0 } t(x)
\in F$ and  $\lim_{t\to \infty}t(x)\in F$.

\begin{definition}[\cW] A quasi-elementary cobordism $B$ is said to be {\em
elementary} if the fixed point set $B^{K^*}$ is connected.
\end{definition}

\begin{definition}[cf. \cite{Morelli1},\cW]
We say that a  birational cobordism $B=B_\phi(X_1,X_2)$ is {\em
collapsible} if the relation $\prec$ is a strict pre-order, namely, there
is no cyclic chain of fixed point components
\[ F_1\prec F_2\prec\ldots\prec F_m\prec F_1.\]
\end{definition}

\subsection{The main example}\label{Subsec:main-example}

We now recall a  fundamental example of an elementary  birational
cobordism in the toric setting, discussed in  \cite{Wlodarczyk2}:

\begin{example}
 Let $B = {\mathbb A}^n =
\Spec K[z_1, \cdots, z_n]$ and let $t\in K^*$ act by
$$t(z_1, \ldots, z_i, \ldots, z_n) = (t^{\alpha_1}z_1,
\ldots, t^{\alpha_i}z_i, \ldots, t^{\alpha_n}z_n).$$ We assume $K^*$ acts
effectively, namely $\gcd(\alpha_1,\ldots,\alpha_n) = 1$.
We regard ${\mathbb A}^n$ as a toric variety defined by a lattice
$N\cong\ZZ^n$ and a regular cone $\sigma\in N_\RR$ generated by the
standard basis
$$\sigma = \langle v_1, \ldots, v_n\rangle.$$
The dual cone $\sigma^\vee$ is generated by the dual basis
$v_1^*,\ldots,v_n^*$, and we identify $z^{v_i^*}=z_i$. 
The $K^*$-action then corresponds to a one-parameter subgroup
$$a = (\alpha_1, \ldots, \alpha_n) \in N.$$
We assume that $a\notin\pm\sigma$. 
We have the obvious description of the sets $B_+$ and $B_-$:
$$\begin{array}{rcl}
B_- &=& \{(z_1, \cdots, z_n); z_i \neq 0 \mbox{ for some }i
\mbox{ with }\alpha_i = (v_i^*,a) < 0\}, \\
B_+ &=& \{(z_1, \cdots, z_n); z_i \neq 0 \mbox{ for some }i
\mbox{ with }\alpha_i = (v_i^*,a) > 0\}.
\end{array}$$
We define the upper boundary and lower boundary fans of $\sigma$ to be 
$$\begin{array}{rcl}
\partial_-\sigma &=& \{x \in \sigma;x + \epsilon \cdot a \not\in \sigma
\mbox{ for all }\epsilon > 0\}, \\ 
\partial_+\sigma &=& \{x \in \sigma;x + \epsilon \cdot (- a) \not\in \sigma
\mbox{ for all }\epsilon > 0\}.
\end{array}$$
Then we obtain the description of $B_+$ and $B_-$ as the toric varieties
corresponding to the fans $\partial_+\sigma$ and $\partial_-\sigma$ in
$N_\RR$.  

Let $\pi: N_\RR \to N_\RR/\RR\cdot a$ be the projection. Then $B/\!/K^*$ is
again an affine toric variety defined by the lattice $\pi(N)$ and
{\em cone} $\pi(\sigma)$. Similarly, one can check that the geometric 
quotients
$B_-/K^*$ and $B_+/K^*$ are toric varieties defined by {\em fans}
$\pi(\partial_+\sigma)$ and $\pi(\partial_-\sigma)$. Since both
$\pi(\partial_+\sigma)$ and $\pi(\partial_-\sigma)$ are subdivisions of
$\pi(\sigma)$, we get a diagram of birational toric maps
$$\begin{array}{ccccc}
B_-/K^* & & \stackrel{\varphi}{\dashrightarrow} & & B_+/K^* \\
& \searrow && \swarrow & \\
&& B /\!/ K^* &&
\end{array} $$

More generally, one can prove that if $\Sigma$ is a subdivision of a convex
polyhedral cone in $N_\RR$ 
with lower boundary $\partial_-\Sigma$ and upper boundary
$\partial_+\Sigma$ relative to an element $a \in
N \setmin \pm \Sigma$, then the toric variety corresponding to
$\Sigma$, with the $K^*$-action given by the one-parameter subgroup $a
\in N$, is a birational cobordism between the two toric varieties
corresponding to 
$\pi(\partial_-\Sigma)$ and $\pi(\partial_+\Sigma)$ as fans in $N_{\mathbb
R}/\RR\cdot a$.

For the details, we refer the reader to \cite{Morelli1},
\cite{Wlodarczyk2} and \cite{Abramovich-Matsuki-Rashid}. \qed
\end{example}

\subsection{Construction of a cobordism} It was shown in \cW\ that birational
cobordisms 
exist for any birational map $X_1 \dashrightarrow X_2$. Here we deal with a
very special case.

\begin{theorem} Let $\phi:X_1 \to X_2$ be a projective birational
morphism between complete nonsingular algebraic varieties, which is
an isomorphism on an open set $U$. Then there is a
complete nonsingular 
algebraic variety $\overline{B}$ with an effective $K^*$-action, satisfying the
following properties:
\begin{enumerate}
\item \label{X-embedded} There exist closed embeddings $\iota_1:X_1
\hookrightarrow
\overline{B}^{K^*}$ 
and  $\iota_2:X_2 \hookrightarrow \overline{B}^{K^*}$ with disjoint images.
\item The open subvariety $B = \overline{B} \setmin \left(\iota_1(X_1)\cup
\iota_2(X_2)\right)$ is a birational cobordism between $X_1$ and $X_2$
respecting the open set $U$.
%\item \label{B-proj}
% There is a $K^*$-equivariant, projective morphism $\overline{B}
% \to X_2 \times \PP^1$, where the action on $X_2 \times \PP^1$ is
% induced by the standard embedding  $K^*\subset \PP^1$.
\item There is a coherent sheaf $E$  on $X_2$, with a $K^*$-action, and a
closed $K^*$-equivariant embedding 
$\overline{B}\subset \PP(E):= {\cP}roj_{X_2}\ \Sym E$.  
\end{enumerate}
\end{theorem}

{\bf Proof.} Let $J\subset \cO_{X_2}$ be an ideal sheaf such that
 $\phi:X_1\to X_2$ is the blowing up morphism 
of $X_2$ along $J$ and
 $J_U=\cO_U$. Let $I_0$ be the ideal of the point
$0\in \PP^1$. Consider $W_0 = X_2 \times \PP^1$ and let $p:W_0
\to X_2 $ and $q:W_0 \to \PP^1$ be the projections. Let $I = (p^{-1} J
+q^{-1} I_0) \cO_{W_0}$. Let $W$ be the blowing up of $W_0$ along $I$.  (Paolo
Aluffi has pointed out that this $W$ is used when constructing the deformation
to the normal cone of $J$.)

We claim that $X_1$ and $X_2$ lie in the nonsingular locus of $W$. For
$X_2\cong X_2\times\{\infty\}\subset X_2\times \bbA^1\subset W$ this is
clear. Since 
$X_1$ is nonsingular, embedded in $W$ as the strict transform of
$X_2\times\{0\}\subset X_2\times \PP^1$, to prove that $X_1$ lies in the
nonsingular locus, it suffices to prove that $X_1$ is a Cartier divisor in
$W$. We look at local coordinates. Let $A=\Gamma(V,\cO_V)$ for some affine
open subset $V\subset X_2$, and let $y_1,\ldots,y_m$ be a set of generators
of $J$ on $V$. Then on the affine open subset $V\times\bbA^1\subset
X_2\times\PP^1$ with coordinate ring $A[x]$, the ideal $I$ is generated by
$y_1,\ldots,y_m,x$. The charts of the blowing up containing the strict
transform 
of $\{x=0\}$ are of the form 
\[ \Spec
A \left[\frac{y_1}{y_i},\ldots,\frac{y_m}{y_i},\frac{x}{y_i}\right] = \Spec
A \left[\frac{y_1}{y_i},\ldots,\frac{y_m}{y_i}\right]\times \Spec
K\left[\frac{x}{y_i}\right], \]  where $K^*$ acts on the second factor.
The strict transform of $\{x=0\}$ is defined by $\frac{x}{y_i}$, hence it is
Cartier. 

Let $\overline{B} \to W$ be a canonical resolution of singularities. Then
conditions 1 and 2 are clearly satisfied. For condition 3, note that
$\overline{B} \to 
X_2\times \PP^1$, being a composition of blowings up of invariant ideals,
admits an 
equivariant ample line bundle. Twisting by the pullback of $\cO_{\PP^1}(n)$ we
obtain an equivariant line bundle which is ample for $\overline{B} \to X_2$.
Replacing this by a sufficiently high power and pushing forward we get $E$.
 \qed 

We refer the reader to \cite{Wlodarczyk2} for more details.

We call a variety $\overline{B}$ as in the theorem a {\em compactified,
relatively projective cobordism}.  

\subsection{Collapsibility and Projectivity} 

Let $B=B_\phi(X_1,X_2)$ be a birational cobordism. We seek a criterion
for collapsibility of $B$.

Let $\cC$ be the set of connected components of
$B_\phi(X_1,X_2)^{K^*}$, and let $\chi:\cC \to \ZZ$ be a function. We
say that $\chi$ is strictly increasing if $F \prec F' \Rightarrow \chi(F)
< \chi(F')$. 
The following lemma is obvious:

\begin{lemma}\label{Lem:order-function}
Assume there exists a strictly increasing function $\chi$. Then 
$\prec$ is a strict  pre-order, and  $B$ is collapsible. Conversely,
suppose $B$ is collapsible. Then there exists a strictly increasing
function $\chi$. \qed
\end{lemma}

\begin{remark} It is evident that every strictly increasing function
can be replaced by one which induces a strict total order. However, it will
be convenient for us to consider arbitrary strictly increasing
functions.
\end{remark}

Let $\chi$ be a strictly increasing function, and let
$a_0<a_1\cdots <a_m\in \ZZ$ be the values of $\chi$. 
\begin{definition} \label{Def:notation-qe} We denote
\begin{enumerate}
\item $F_{a_i}=\cup\{F|\chi(F)=a_i\}$.
\item  $F_{a_i}^+=\cup\{F^+|\chi(F)=a_i\}$.
\item  $F_{a_i}^-=\cup\{F^-|\chi(F)=a_i\}$.
\item  $F_{a_i}^\pm=\cup\{F^\pm|\chi(F)=a_i\}$.
\item  $F_{a_i}^*=\cup\{F^*|\chi(F)=a_i\}$.
\item  $B_{a_i}=B \setmin
\left(\ \bigcup\{F^-|\chi(F)<a_i\}\quad \cup
\quad\bigcup\{F^+|\chi(F)>a_i\}\ \right)$. 
\end{enumerate}
\end{definition}

Note that $F_{a_i}^*$ is the union of non-closed orbits in $B_{a_i}$. The
following is an immediate extension of Proposition 1 of \cW.   

\begin{proposition}
\begin{enumerate}
\item $B_{a_i}$ is a quasi-elementary cobordism.
\item For $i=0,\ldots, m-1$ we have $(B_{a_i})_+ = (B_{a_{i+1}})_-$.
\end{enumerate}
\end{proposition} 

The following is an analogue of Lemma 1 of \cW\ in the
case of the cobordisms we have constructed.

\begin{proposition} Let 
$E$ be a coherent sheaf on $X_2$ with a $K^*$-action, and let
$\overline{B} \subset {\mathbb P}(E)$ be a  
compactified, relatively projective cobordism  embedded $K^*$-equivariantly.
Then there exists a strictly increasing function $\chi$ for the cobordism
$B = \overline{B} \setmin(X_1 \cup X_2)$. In particular, the  cobordism is
collapsible.
\end{proposition}

{\bf Proof.} Since $K^*$ acts trivially on $X_2$, and since $K^*$ is reductive,
there exists a direct sum decomposition 
$$ E = \mathop{\bigoplus}\limits_{b\in \ZZ} E_b$$
where $E_b$ is the subsheaf on which the action of $K^*$ is given by
the character $t\mapsto t^b$. 
%(Proof: If $K^d\times V\to V$ is a
%trivialization of $E$ over an invariant affine open subset $V\subset
%X_2$ then the action of $K^*$ on $E$ in this trivialization is given
%by a group homomorphism $K^*\to \Gamma(V,GL_d(K))$. Since $K^*$ is
%reductive, we can diagonalize all these matrices simultaneously.) 
Denote by $b_0,\ldots,b_k$ the characters which figure in
this representation.
%, and let $d_i=\operatorname{rank} E_{b_i}$.
Note that there are disjoint embeddings $\PP(E_{b_j}) \subset \PP(E)$. 

Let $p\in B$ be a fixed point lying in the fiber $\PP(E_q)$ over $q\in X_2$.
We choose a basis 
$$(x_{b_0,1},\ldots, x_{b_0,d_0},\ldots,x_{b_k,1},\ldots, x_{b_k,d_k})$$
 of $E_q$ where $x_{b_j,\nu}\in E_{b_j}$ and 
use the following lemma:
\begin{lemma}
Suppose $p\in \PP(E_q)^{K^*}$ is a fixed point with homogeneous  coordinates
$$(p_{b_0,1},\ldots, p_{b_0,d_0},\ldots,
p_{b_k,1},\ldots, p_{b_k,d_k}).$$ Then  there is an $j_p$ such that
$p_{b_j,\nu}=0$ whenever $j\neq j_p$. In particular, $p\in
\PP(E_{b_{j_p}})\subset \PP(E)$. \qed
\end{lemma}

If $F\subset B^{K^*}$ is a connected component of the fixed point set,
then it follows from the lemma that  $F \subset \PP(E_{b_j})$ for some
$j$. We define $$\chi(F) = b_j.$$
To check that $\chi$ is strictly increasing, consider a point $p\in B$ such
that $\lim_{t\to 0} t(p) \in F_1$ and $\lim_{t\to \infty} t (p) \in F_2$ for
some fixed point components $F_1$ and $F_2$.  
Let the coordinates of $p$ in the fiber over $q\in X_2$ be $(p_{b_0,1},\ldots,
p_{b_0,d_0},\ldots, p_{b_k,1},\ldots, p_{b_k,d_k})$. Now 
\begin{align*}
\lim_{t\to 0} t(p) &\in \PP(E_{b_{min}}),\\
\lim_{t\to \infty} t(p) &\in \PP(E_{b_{max}}),
\end{align*}
where 
\begin{align*}
b_{min}=  \min\{b_j: p_{b_j,\nu}\neq 0 \mbox{ for some } \nu\},\\
b_{max}=  \max\{b_j: p_{b_j,\nu}\neq 0 \mbox{ for some } \nu\}.
\end{align*}
Thus, if $p$ is not fixed by $K^*$ then
\[ \chi(F_1) = b_{min} < b_{max} = \chi(F_2).\]
\qed

\subsection{Geometric invariant theory and projectivity}
In this section we use geometric invariant theory, and ideas
(originating in symplectic geometry) developed by
M. Thaddeus and others (see
e.g. \cite{Thaddeus2}), in 
order to obtain a result about relative projectivity of quotients.

We continue with the notation of the last section. Consider the sheaf $E$ and
its decomposition according to the character. Let $\{b_j\}$ be 
the characters of the
action of $K^*$ on $E$, and $\{a_i\}$ the subset of those
$b_j$ that are in the image of $\chi$.
If we use the Veronese embedding $\overline{B} \subset \PP(\Sym^2(E))$ and
replace $E$ by $\Sym^2(E)$, we may assume that $a_i$ are even, in
particular $a_{i+1} > a_i+1$ (this is a technical condition which
comes up handy in what follows). 

Denote by $\rho_0(t)$ the action of $t\in K^*$ on $E$. For any $r\in \ZZ$
consider the ``twisted'' action $\rho_r(t)=
t^{-r}\cdot\rho_0(t)$. Note that the induced action on $\PP(E)$ does
not depend on the ``twist'' $r$. Considering the decomposition
$E = \mathop{\bigoplus} E_{b_j}$, we see that $\rho_r(t)$ acts
on $E_{b_j}$ by multiplication by $t^{b_j-r}$. 

We can apply geometric invariant theory in its relative form (see,
e.g., \cite{Pandharipande}, \cite{Hu2}) to the action 
$\rho_r(t)$ of $K^*$. Recall that a point $p\in \PP(E)$ is said to be
semistable with respect to $\rho_r$, written $p\in
(\PP(E),\rho_r)^{ss}$, if there is a positive integer $n$ and a
$\rho_r$-invariant local section
$s\in (\Sym^n(E))^{\rho_r}$, such that $s(p)\neq 0$.
The main result of geometric invariant theory implies that 
$$\mathop{{\cP}roj}\limits_{X_2}\
\mathop{\bigoplus}\limits_{n\geq 0}^\infty 
(\Sym^n(E))^{\rho_r} = (\PP(E),\rho_r)^{ss} /\!/ K^*;$$
moreover, the quotient map $(\PP(E),\rho_r)^{ss} \to (\PP(E),\rho_r)^{ss} /\!/
K^*$ is affine. 
We can define $(\overline{B},\rho_r)^{ss}$ analogously, and we automatically
have $(\overline{B},\rho_r)^{ss} = \overline{B} \cap(\PP(E),\rho_r)^{ss}$.

The numerical criterion of semistability (see
\cite{Mumford-Fogarty-Kirwan}) immediately implies the following:
\begin{lemma} For $0<i<m$ we have  
\begin{enumerate}
\item $(\overline{B},\rho_{a_i})^{ss} = B_{a_i}$;
\item $(\overline{B},\rho_{a_i+1})^{ss} = (B_{a_i})_+$;
\item $(\overline{B},\rho_{a_i-1})^{ss} = (B_{a_i})_-$;
\end{enumerate}
\end{lemma}

In other words, the triangle of birational maps 
$$ \begin{array}{ccccc}
 ({B}_{a_i})_-/K^* & & \stackrel{\varphi_i}{\das} & &({B}_{a_i})_+/K^* \\
 & \searrow & & \swarrow & \\
& & {B}_{a_i} /\!/ K^* & & \end{array}$$
is induced by by a change of linearization of the action of $K^*$.

In particular we obtain:

\begin{proposition}\label{prop:projectivity}
The morphisms $({B}_{a_i})_+/K^*\to X_2$, $({B}_{a_i})_-/K^*\to X_2 $ and  
${B}_{a_i} /\!/ K^*\to X_2 $ are projective.
\end{proposition}

\subsection{The main result of \cite{Wlodarczyk2}}

Let $B$ be a collapsible nonsingular birational
cobordism. Then we can write $B$ as a union of quasi-elementary cobordisms
$B=\cup_i B_{a_i}$, with $(B_{a_i})_+ = (B_{a_{i+1}})_-$. By
Lemma~\ref{lem:loc-tor-act} each
$B_{a_i}$ has a locally toric structure such that the action of $K^*$
is locally toric.

\begin{lemma}\label{lem:loc-descr} Let $B_{a_i}$ be a 
quasi-elementary cobordism, with a relatively affine {\em locally toric} $K^*$ action. Then
$B_{a_i}/\!/K^*$, $(B_{a_i})_-/K^*$, 
$(B_{a_i})_+/K^*$ are locally toric varieties and we have a diagram of
locally toric maps 
$$\begin{array}{ccccc}
(B_{a_i})_-/K^* & & \stackrel{\varphi_i}{\dashrightarrow} & & (B_{a_i})_+/K^*
\\
& \searrow && \swarrow & \\
&& B_{a_i} /\!/ K^* &&
\end{array} $$
 where $\varphi_i$ is a tightly
locally toric birational map.

In case $B_{a_i}$ is nonsingular, the diagram above can be described in toric
charts by the main example in  Section~\ref{Subsec:main-example}.  

If the action of $K^*$ on $B_{a_i}$ is toroidal then all these varieties and
maps are also toroidal, and  $\varphi_i$ is a tightly toroidal birational map.
\end{lemma}

{\bf Proof.} Let $\eta_p:V_p\to X_p$ be a strongly \'etale $K^*$-equivariant
toric chart in $B_{a_i}$ giving a locally toric structure to the action of 
$K^*$. Then $(V_p)_- = (B_{a_i})_-\cap V_p$ and the morphism $(V_p)_- \to
(X_p)_-$ is again strongly \'etale, providing locally toric structures on the
variety $(B_{a_i})_-/K^*$ and the morphism  $(B_{a_i})_-/K^* \to
B_{a_i}/\!/K^*$. Similarly for $(B_{a_i})_+$. \qed 

Now we assume $B\subset \overline{B}$ is open in a compactified, relatively
projective cobordism.
When we compose the birational transformations obtained from each $B_{a_i}$
we get a slight refinement of the main result of \cW.

\begin{theorem}
\label{Th:locally-toric-factorization}
 Let
 $\phi:X_1 \dashrightarrow X_2$ be a birational map between complete
 nonsingular algebraic varieties $X_1$ and $X_2$ over an algebraically closed
 field $K$ of characteristic zero, and let $U\subset  X_1$ be an
 open set where $\phi$ is an isomorphism.  Then there exists a 
 sequence of birational maps between complete locally toric algebraic 
varieties  
 $$X_1 = W_0 \stackrel{\varphi_1}{\dashrightarrow} W_1
 \stackrel{\varphi_2}{\dashrightarrow} \cdots
 \stackrel{\varphi_i}{\dashrightarrow} W_i
 \stackrel{\varphi_{i+1}}{\dashrightarrow} W_{i+1}
 \stackrel{\varphi_{i+2}}{\dashrightarrow}
 \cdots \stackrel{\varphi_{m-1}}{\dashrightarrow}
 W_{m-1} \stackrel{\varphi_m}{\dashrightarrow} W_m = X_2$$   where
 \begin{enumerate}
   \item  
   $\phi = \varphi_m \circ \varphi_{m-1} \circ \cdots \varphi_2 \circ
 \varphi_1$,
   \item $\varphi_i$ are isomorphisms on $U$, and 
   \item For each $i$, the map $\varphi_i$ is tightly locally toric, and 
     \'etale  locally equivalent to a map $\varphi$ described in
     \ref{Subsec:main-example}. 
 \end{enumerate}
Furthermore, there is an index $i_0$ such that for all $i\leq i_0$ the map
$W_i\das X_1$ is a projective morphism, and for all $i\geq i_0$ the map
$W_i\das X_2$ is a projective morphism. In particular, if $X_1$ and
$X_2$ are  projective  then all the $W_i$ are  projective.
\end{theorem}

\begin{remark} For the projectivity claim 2, we take the first $i_0$ terms
in the factorization to come from Hironaka's elimination of
indeterminacies in Lemma~\ref{lem:red-proj}, which is projective over $X_1$,
whereas the last 
terms come from  $\overline{B}$, which is projective over $X_2$, and
the geometric invariant theory considerations as in
Proposition~\ref{prop:projectivity}.  
\end{remark} 

\subsection{Projectivity of toroidal weak
factorization}\label{Sec:proj-tor-fac} 
The following is a refinement of Theorem \ref{thm:morelli}, in which a
projectivity statement is added:

\begin{theorem}\label{thm:morelli-proj} Let $U_1\subset W_1$ and
$U_2\subset W_2$ be nonsingular toroidal 
embeddings.
Let $\psi: W_1\dra W_2$ be a proper toroidal birational
map. Then $\phi$ can be factored into a
sequence of toroidal birational maps consisting of smooth toroidal  blowings
up and down, namely:
 $$
 W_1 = V_0 \stackrel{\varphi_1}{\dashrightarrow} V_1
 \stackrel{\varphi_2}{\dashrightarrow} \cdots
 \stackrel{\varphi_i}{\dashrightarrow} V_i
 \stackrel{\varphi_{i+1}}{\dashrightarrow} V_{i+1}
 \stackrel{\varphi_{i+2}}{\dashrightarrow}
 \cdots \stackrel{\varphi_{l-1}}{\dashrightarrow}
 V_{l-1} \stackrel{\varphi_l}{\dashrightarrow} V_l = W_2
$$
   where
 \begin{enumerate}
   \item  
   $\phi = \varphi_l \circ \varphi_{l-1} \circ \cdots \varphi_2 \circ
 \varphi_1$;
   \item $\varphi_i$ are isomorphisms on $U$, the embeddings  $U \subset V_i$
 are toroidal, and $\varphi_i$ are toroidal birational maps; and
  \item either 
  $\varphi_i:V_i \dashrightarrow V_{i+1}$ or $\varphi_i^{-1}:V_{i+1}
\dashrightarrow V_{i}$ 
  is a toroidal morphism obtained by blowing up a smooth irreducible toroidal
 center. 
 \end{enumerate}
Furthermore, there is an index $i_0$ such that for all $i\leq i_0$ the map
$V_i\das X_1$ is a projective morphism, and for all $i\geq i_0$ the map
$V_i\das X_2$ is a projective morphism. In particular, if $X_1$ and
$X_2$ are  projective  then all the $V_i$ are  projective. 
\end{theorem}

{\bf Proof.} As in \cite{Abramovich-Matsuki-Rashid}, Lemma 8.7
 we reduce to the case where the
polyhedral complex of $W_2$ is embeddable as a quasi-projective toric fan
$\Delta_2$ in a space $N_\RR$.  
Indeed that Lemma  gives an embedding preserving the $\QQ$-structure for the
barycentric subdivision of any simplicial complex, and since $\Delta_2$ is
nonsingular this embedding preserves integral structures as
well. A further subdivision ensures that the fan is quasi-projective. (We note
 that this embedding is introduced for the sole purpose of 
applying Morelli's $\pi$-Desingularization Lemma directly, rather than 
observing  that the proof works word for word in the toroidal case.)

  As in \ref{thm:morelli} we may assume $W_1 \das W_2$ is a
projective morphism. Thus the complex $\Delta_1$ of $W_1$ is a projective
 subdivision of $\Delta_2$. 
Our 
construction of a compactified relatively projective cobordism $\overline B$
for the morphism $\phi$ yields a toroidal embedding $B$ whose complex
$\Delta_B$ is a quasi-projective polyhedral cobordism lying in  $(N\oplus
\ZZ)_\RR$ such that $\pi(\partial_+\Delta_B) = \Delta_2$ and
$\pi(\partial_-\Delta_B) = \Delta_1$, where $\pi$ is the projection onto
$N_\RR$. Moreover, the toroidal morphism 
$B \to W_2$ gives a polyhedral morphism $\Delta_B \to \Delta_2$ induced by the
projection $\pi$. Morelli's
$\pi$-desingularization lemma gives a projective subdivision $\Delta_B'
\to\Delta_B$, isomorphic on the upper and lower boundaries
$\partial_\pm\Delta_B$, such that $\Delta_B'$ is $\pi$-nonsingular. We still
have 
a polyhedral morphism $\Delta_B' \to \Delta_2$. The complex $\Delta_B'$
corresponds to a toroidal birational cobordism $B'$ between $W_1$ and
$W_2$. Since $\Delta_B'$ is $\pi$-nonsingular, any elementary piece
$B_F'\subset 
B'$ corresponds  a toroidal blowing up followed by a toroidal blowing down
between 
nonsingular toroidal embeddings, with 
nonsingular centers. It follows that the same holds for every
quasi-elementary piece of $B'$ (here the centers may be reducible). As in 
Theorem 
\ref{Th:locally-toric-factorization} above, these toroidal embeddings can be
chosen to be
projective over $W_2$. \qed

\section{Torification}\label{Sec:torification}

We wish to replace the locally toric factorization of Theorem 
\ref{Th:locally-toric-factorization} by a toroidal
factorization. This amounts to replacing $B$ with a locally toric
$K^*$-action by some $B'$ with a toroidal $K^*$-action. We call such a
procedure {\em torification}. The basic idea, which goes back at least to 
Hironaka, is that if 
one blows up an ideal, the exceptional divisors provide the resulting
variety with useful extra structure. The ideal we construct, called a
{\em torific ideal},  is closely related to the torific ideal of
\cite{Abramovich-de-Jong}. 

\subsection{Construction of a torific ideal}\label{Sec:torific-construction}
\begin{definition} Let $V$ be an 
algebraic variety with a $K^*$-action, $p\in V$ a closed point,
$G_p\subset K^*$ the stabilizer of $p$. Fix an integer $\alpha$.
Then we define
$$ J_{\alpha,p} \subset {\cO}_{V,p}$$
to be the ideal generated by the semi-invariant functions $f \in
{\cO}_{V,p}$ of $G_p$-character $\alpha$, that is, for $t\in G_p$
%and a generator $f\in J_{\alpha,p}$
 we have
$$t^*(f) = t^{\alpha}f.$$
\end{definition}

\begin{lemma}\label{lem:ideal-desc} Let $V$ be a variety with a
$K^*$-action, and $p\in V$ a closed point. If $z_1,\ldots,z_n$ are
$G_p$-semi-invariant generators of the maximal ideal  ${\mathfrak m}_p$, then
$J_{\alpha,p}$ is generated by monomials in $z_i$ having $G_p$-character
$\alpha$. 
\end{lemma}

{\bf Proof.} Consider the completion of the local ring $\hat{\cO}_{V,p}$, 
and lift the action of $G_p$ to it.
Since  $\hat{\cO}_{V,p}$ is a faithfully flat $O_{V,p}$-module,
it suffices to prove  
that the completion $\hat{J}_{\alpha,p}$ is the ideal of $\hat{\cO}_{V,p}$
generated by monomials in $z_i$ of $G_p$-character $\alpha$.

Consider the $G_p$-equivariant epimorphism $$K[\![z_1,\ldots,z_n]\!] \to
\hat{\cO}_{V,p}.$$
Since $G_p$ is reductive,  a semi-invariant element of $\hat{\cO}_{V,p}$ is the
image of a semi-invariant power series in $z_i$.
A monomial in $K[\![z_1,\ldots,z_n]\!]$ of $G_p$-character $\alpha$ clearly
maps to 
$\hat{J}_{\alpha,p}$. Conversely, a semi-invariant power series in $z_i$ must
have all its monomials semi-invariant of the same character.
One can choose a finite set of monomials
occurring in the power series such that any other monomial occurring in this
power series is divisible by one of them. Hence the power series lies in the
ideal generated by monomials in $z_i$ of $G_p$-character $\alpha$. \qed

The lemma implies that, given a strongly \'etale morphism $\eta: V\to X$
between varieties with $K^*$ action, the
inverse image of $J_{\alpha,\eta(p)}$ generates $J_{\alpha,p}$. Indeed, we can
choose $G_p=G_{\eta(p)}$ semi-invariant generators of  ${\mathfrak m}_\eta(p)$,
which pull 
back to semi-invariant generators of ${\mathfrak m}_p$.

For the rest of this section, we let $B$ be a {\em 
quasi-elementary} cobordism with a {\em relatively affine, locally toric
$K^*$-action}; $B=B_{a_i}$ for some $i$ according to our 
previous notation. Without loss of generality we assume $a_i = 0$,
so $F_0 = B^{K^*}$. Recall the notation $F_0^*$ in Definition
\ref{Def:notation-qe}. This is a 
constructible set in $B$, which is the union of the non-closed orbits.

\begin{proposition}\label{Prop:ideal-existence} There exists a unique
coherent $K^*$-equivariant ideal sheaf 
$I_\alpha$ on $B$, such that for all $p\in B\setmin F_{0}^*$ we have
$(I_\alpha)_p = J_{\alpha,p}$. 
\end{proposition} 

\begin{definition} 
The sheaf $I_\alpha$ is called {\em the $\alpha$-torific ideal sheaf} of the
action of $K^*$ on $B$.
\end{definition}

\begin{remarks} 
\begin{enumerate}
\item Notice that the collection of ideals $J_{\alpha,p}$ for
$p\in B$ does not define a coherent sheaf of ideals in general.
As an example, let $B=\bbA^2$, and let  $t\in K^*$ act by 
\[ t(x,y)= (t x,t^{-1} y).\]
Then at $p=(0,0)$, the stabilizer is $G_p = K^*$, and $J_{1,p}=(x)$. Any other
point $q\in B\setmin\{p\}$ has a trivial stabilizer, hence $J_{1,q} =
\cO_{B,q}$ is 
trivial. These germs do not form a coherent ideal sheaf on
$B$. In this case, the ideal sheaf generated by $x$ is the $1$-torific ideal
sheaf 
$I_1$ of the proposition. 
\item Note also that the assertion of the proposition fails if we remove the
requirement on $B$ being quasi-elementary. For a simple example which is not a
coobrdism, let  $t\in K^*$ act on
$B = \PP^1$ with homogeneous coordinates $(X:Y)$, via $(X:Y) \to (tX:Y)$. Then
at $p=(1:0)$ the ideal $J_{1,p}$ is generated by $x_0 = X/Y$, whereas
at 
$q=(0:1)$ the ideal $J_{1,q}$ is the zero ideal. It is easy to construct
higher 
dimensional examples of cobordisms where the ideal $J_{\alpha,p}$ at one fixed
point cannot 
be 
glued  to {\em any} $\beta$-torific ideal at another fixed point.
\end{enumerate}
\end{remarks}

{\bf Proof of \ref{Prop:ideal-existence}.}
First we prove the uniqueness of $I_\alpha$. Clearly the stalks of $I_\alpha$
are uniquely defined at all points $p\in B\setmin F_{0}^*$. Now if $p\in
F_0^*$ then $p$ has a unique limit fixed point  $p'\in F_0$, where either
$p'=\lim_{t\to 0}t(p)$, or $p'=\lim_{t\to \infty}t(p)$, but not both, since $B$
is quasi-elementary. Since
$I_\alpha$ is 
uniquely determined at $p'$, hence also near $p'$ by coherence, it follows
from $K^*$-equivariance that $I_\alpha$ is uniquely determined at $p$.

To prove the existence of $I_\alpha$ we cover $B$ with strongly \'etale affine
toric 
charts $\eta_p:V_p\to X_p$ as in Lemma~\ref{lem:loc-tor-act}. With such
charts, it follows that 
$F_0^*$ restricted to $V_p$ is the inverse image of $F_0^*$ defined in $X_p$
(recall that $F_0^*$ consists of the union of non-closed orbits).
Lemma \ref{Lem:ideal-existance-toric} below gives the existence
of a torific ideal on $X_p$, and its pullback is a torific ideal on $V_p$  by
Lemma~\ref{lem:ideal-desc}. By uniqueness, the ideals defined on $V_p$ glue
together to an ideal on $B$. \qed

\begin{lemma}\label{Lem:ideal-existance-toric} Let $B=X(N,\sigma)$ be an
affine toric variety on which $K^*$ acts as a one-parameter subgroup of the
torus. 
Then $I_\alpha$ exists and is generated by all monomials $z^m, m\in\sigma^\vee$
on which $K^*$ acts
by character $\alpha$. 
\end{lemma}

{\bf Proof.} 
By abuse of notation we will use the same letter
$\alpha$ to denote a character of a subgroup of $K^*$. 
It follows easily from Lemma \ref{lem:ideal-desc} that
 for any $p\in B$ the ideal $J_{\alpha,p}$ is
generated by all elements $z^m$ regular at $p$ on which $G_p$ acts by
character $\alpha$. 

Let $p\in B\setmin F_0^*$, and let $\tau$ be the
smallest face of $\sigma$ such that $p$ lies in the affine open toric
subvariety $X(N,\tau)$. Then the  monomials $z^m$ regular at $p$ are
those for which $m\in M\cap\tau^\vee$, and the monomials invertible at $p$
are the 
ones for which  $m\in M\cap\tau^\perp$. 

If $z^m$ for $m\in M\cap\sigma^\vee$ is a monomial regular on B on which
$K^*$ acts 
by character $\alpha$ then clearly $z^m$ is regular at $p$ and
$G_p\subset K^*$ acts on it by character $\alpha$. 
Conversely, let $z^m$ for $m\in M\cap\tau^\vee$ be a monomial regular at $p$,
on which $G_p$ acts by character $\alpha$. We show that there exists a
monomial $z^{m'}$ invertible at $p$ 
(i.e., $m'\in M\cap \tau^\perp$) such that $z^{m+m'}$ is
regular on $B$ and $K^*$ acts on it by character $\alpha$. This is done in
two steps: 

STEP 1. There exists $m'\in M\cap \tau^\perp$ such that $z^{m+m'}$ has
$K^*$-character $\alpha$. Since $G_p$ is the subgroup of $K^*$ acting
trivially on the monomials corresponding to $m'\in M\cap\tau^\perp$, we have
an exact sequence  
\[ M\cap\tau^\perp \to \hat{K}^* \to \hat{G}_p \to 0, \]
where $\hat{H}$ denotes the character group of $H$. Thus, we may replace $m$
by $m+m'$ and assume that $K^*$ acts on $z^m$ by character $\alpha$. 

STEP 2. There exists $m'\in M\cap\tau^\perp$ such that $z^{m+m'}$ is regular
on $B$, i.e., $m+m'\in M\cap
\sigma^\vee$. Since the monomial $z^m$ is $K^*$-semi-invariant, there exists
an 
affine open $K^*$-invariant neighborhood of $p$ on which $z^m$ is regular,
and 
since $p\notin F_0^*$, this neighborhood can be chosen of the form
$\pi^{-1}(V)$ where $\pi: B\to B/\!/K^*$ is the projection and $V\subset
B/\!/K^*$ is an affine open toric subvariety. Let $m'\subset M\cap
\sigma^\vee\cap 
a^\perp$ be such that $V$ is the nonvanishing locus
of the monomial $z^{m'}$. Then $z^{m'}$ as a monomial on $B$ is invertible at
$p$, has 
$K^*$-character $0$, and replacing $m'$ by a multiple if necessary, we have
that $z^{m+m'}$ is regular on $B$. \qed

\subsection{The torifying property of the torific
ideal}\label{Sec:torifying-property} 

Suppose $X$ is a locally toric variety with a locally toric action of
$K^*$, and $D\subset X$ a divisor compatible with the locally toric
structure, that means, at each point of $X$ we can find a toric chart
$\eta_p:V_p\to X_p$ such that $D\cap V_p$ is the inverse image of some
toric divisor $D_p\subset X_p\setmin T$. In this situation we want
to know if $(X \setmin D) \subset X$ is a toroidal embedding  on which
$K^*$ acts 
toroidally. Clearly it suffices to show that $(X_p\setmin D_p) \subset X_p$
 is a toroidal embedding with a toroidal $K^*$ action for all toric charts for
the $K^*$ action. The 
following 
lemma gives several equivalent conditions for this. 

\begin{lemma} \label{Lem:remove-divisors} Let $X=X(\Sigma, N)$ be a toric
variety, $D\subset 
X\setmin T$ a divisor in $X$, and let $K^*$ act on $X$ as a
one-parameter subgroup of the torus $T$, corresponding to a lattice
point $a\in N$. Then the following are equivalent:
\begin{enumerate}

\item $X\setmin D\subset X$ is a toroidal embedding on which $K^*$ acts
toroidally. 

\item For every affine open toric subvariety $X_\sigma\subset X$
corresponding to a cone $\sigma\in \Sigma$, there exists a toric
variety $X_{\sigma'}$ with an action of $K^*$ as a one-parameter
subgroup of the torus $T'$ such that we have a decomposition
\begin{align*}
X_\sigma &\cong  \bfa^k \times  X_{\sigma'} \\
D        &\cong  \bfa^k \times  (X_{\sigma'}\setmin T'),
\end{align*}
where the action of $K^*$ on $X_\sigma$ is a product of the action on
$X_{\sigma'}$ with the trivial action on $\bfa^k$.

\item For every cone $\sigma = \langle v_1,\ldots,v_m\rangle \in
\Sigma$, with $v_1,\ldots, v_{k}$ corresponding to the irreducible
toric divisors not in $D$, we have a decomposition
\begin{align*}
\sigma &\cong \langle v_1, \ldots, v_{k}\rangle \times \langle
v_{k+1},\ldots, v_m\rangle \\
N &\cong N'\times N'',
\end{align*}
where $N'\subset N$ is the sublattice generated by $v_{1},\ldots,
v_k$, and $N''\subset N$ is a complementary sublattice containing $v_{k+1},
\ldots, v_{m}$ as well as the point $a$.

\item \label{It:decomp-cones-single} For every cone $\sigma = \langle
v_1,\ldots,v_m\rangle \in 
\Sigma$, and every $v_i$ corresponding to an irreducible
toric divisor not in $D$, we have a decomposition
\begin{align*}
\sigma &\cong \langle v_i\rangle \times 
\langle v_1, \ldots,\widehat{v}_i,\ldots, v_m\rangle  \\ 
N &\cong N_i\times N_{\widehat\iota},
\end{align*}
where $N_i\subset N$ is the sublattice generated by $v_i$, and $N_{\widehat
\iota}\subset 
N$ is a complementary sublattice containing $v_j$, $j\neq i$ as well as the
point $a$. 

\item For every affine open toric subvariety $X_\sigma\subset X$
corresponding to a cone $\sigma\in \Sigma$, and every irreducible
toric divisor $E$ in $X$ not in $D$, there exists a toric
variety $X_{\sigma'}$ with an action of $K^*$ as a one-parameter
subgroup of the torus $T'$ such that we have a decomposition
\begin{align*}
X_\sigma &\cong \bfa^1 \times X_{\sigma'}  \\
E        &\cong \{0\}  \times X_{\sigma'},
\end{align*}
where the action of $K^*$ on $X_\sigma$ is a product of the action on
$X_{\sigma'}$ with the trivial action on $\bfa^1$.
\end{enumerate}
\end{lemma} 
\begin{remark} We say that a divisor $E$ as in condition 5 is {\em removed from
the toroidal structure of $T\subset X$}.
\end{remark}

{\bf Proof.} The equivalences $2 \Leftrightarrow 3$ and
$4\Leftrightarrow 5$  are simply translations 
between toric varieties and the corresponding fans. The equivalence
$3\Leftrightarrow 4$ follows easily  from the combinatorics of cones. For $2
\Rightarrow 1$, we cover $\bfa^k \times X_{\sigma'}$ with toroidal charts of
the form 
$\bfg_m^k \times X_{\sigma'}$. 
The converse will not be used in this paper, and we leave it to the reader.
%The converse follows from the fact that if
%$X_\sigma$ cannot be decomposed as $X_{\sigma'} \times \bfa^k$ then it is
%not equisingular along the intersection of all divisors in $D$. 
\qed

As before, let $B$ be a   quasi-elementary cobordism, with a relatively affine,
locally toric $K^*$ action. We further assume that $B$ is {\em
nonsingular}. Choose 
$c_1,\ldots,c_\mu\in 
\ZZ$ a finite set of integers representing all characters of the $G_p$-action
on the tangent space of $B$ at $p$ for all $p\in B$. Let 
\[ I= I_{c_1}\cdots I_{c_\mu}\]
be the product of the $c_i$-torific ideals, and let $B^{tor}\to B$ be the
normalized blowing up of $B$ along $I$. Since $I$ is $K^*$-equivariant, we can
lift the action of $K^*$ to $B^{tor}$. Denote by $D\subset B^{tor}$ the total
transform of the support of $I$, and $U_{B^{tor}} = B^{tor}\setmin D$.

We remark that such a set $\{c_1,\ldots, c_\mu\}$ can be found
by covering the quasi-elementary cobordism $B$ with a finite number of
toric charts, and collecting all characters of the $K^*$-action on the
coordinates of the toric varieties appearing in these charts. 
 We are also allowed to enlarge the  set of $c_i$ - this will be utilized in
the following section.

To understand the morphism  $B^{tor} \to B$ we use the following easy lemma:
\begin{lemma}\label{Lem:ideal-nonzero} For any $c\in \ZZ$, the ideal $I_c$ is
nonzero.
\end{lemma}
{\bf Proof.}  This can be seen from the toric
picture given in Lemma~\ref{Lem:ideal-existance-toric}: the lattice point  $a$
is primitive (since the $K^*$-action is effective), 
and 
$\pm a$ do not lie in the cone $\sigma$ (since $B_-\cap B_+$ is a nonempty
open 
set), therefore the hyperplane $(a,\cdot) = c$
contains lattice points in $\sigma^\vee$. Thus the set of $f\in 
{\cO}_{V,p}$ of $G_p$-character $\alpha$ is nonempty. \qed

It follows that $B^{tor}$, being the normalized blowing up of the product
$I_{c_1}\cdots I_{c_\mu}$, satisfies a universal property: it is the minimal
normal modification of $B$ such that the inverse image of $I_{c_i}$ is
principal for all $i$. This implies that $B^{tor}$ is canonically isomorphic
to the normalization of the variety obtained from $B$ by first blowing up
$I_{c_1}$, then the inverse image of $I_{c_2}$, and so on.

\begin{proposition}\label{prop:torifying} The variety $B^{tor}$ is a
quasi-elementary cobordism, with $(B^{tor})_+ = B^{tor} \times_BB_+$ and
$(B^{tor})_- = B^{tor} 
\times_BB_-$. Moreover, the 
embedding  $U_{B^{tor}} \subset B^{tor}$ is  toroidal and $K^*$ acts toroidally
on this embedding. 
\end{proposition}

\begin{definition} We call $I$ a {\em torific} ideal and $B^{tor}\to B$ a
{\em torific} blowing up. 
\end{definition}

{\bf Proof.} Suppose $B^{tor}$ is not quasi-elementary, that means, there
exists a non-constant orbit with both of its limits in $B^{tor}$. Since $B$ is
quasi-elementary and the morphism $B^{tor}\to B$ is equivariant, the image of
this orbit must be a fixed point. However, the coordinate ring of an affine
chart in a
$K^*$-invariant fiber of the morphism $B^{tor}\to B$ is generated by
fractions $f = f_1/f_2$ where $f_i$ are generators of the ideal $I$, hence
$K^*$ acts trivially on $f$. This means that the fiber consists of fixed
points, a contradiction. 

Since $(B^{tor})^{K^*}$ is the inverse image of $B^{K^*}$, we get that $x\in
(B^{tor})_+$ if and only if its image is in $B_+$, and similarly for
$(B^{tor})_-$. 

To prove that $U_{B^{tor}}\subset B^{tor}$ is toroidal and $K^*$ acts
toroidally on this embedding, we consider toric charts $\eta_p: V_p\to X_p$
in $B$ giving the action of $K^*$ on $B$ a locally toric structure. By
Lemma~\ref{lem:ideal-desc} the ideal $I$ restricted to $V_p$ is the inverse
image of the ideal $I_p = I_{p,c_1}\cdots I_{p,c_\mu}$ in $X_p$. It follows
that the normalization of the blowing up of $I_p$ in $X_p$ provides a toric
chart for $B^{tor}$ 
such that the action of $K^*$ on $B^{tor}$ is again locally toric. Let $X_q$ be
such a chart:
$$\begin{array}{ccccc}
X_q & \leftarrow & V_q & \subset & B^{tor} \\
\dar & & \dar & & \dar \\
X_p & \leftarrow & V_p & \subset & B.
\end{array}$$
Let $D_q\subset X_q$ be the support of the divisor defined by the total
transform of $I_p$. Then 
\[ U_{B^{tor}}\cap V_q = \eta_q^{-1} X_q\setmin D_q,\]
and we are reduced to proving that $(X_q\setmin D_q)\subset X_q$ is a
toroidal embedding on which $K^*$ acts toroidally. We do this by verifying the
equivalent condition 5 
in  Lemma \ref{Lem:remove-divisors}.

Let $X_p=X(N,\sigma)$ be a nonsingular affine toric variety defined by the cone
$\sigma 
= \langle v_1,\ldots,v_m\rangle$, 
$\sigma^\vee=\langle v_1^*,\ldots,v_m^*,\pm v_{m+1}^*,\ldots,\pm v_n^*
\rangle$, and let $K^*$ act on $z_i$
by character $c_i$. The only irreducible toric divisors in $X_q$ that do
not lie in
the total transform of $I_p$ are among the strict transforms of the divisors
$\{z_i=z^{v_i^*} = 0\}\subset X_p$. Consider the divisor $\{z_1=0\}$. The 
ideal $I_{p,c_1}$ contains $z_1$. If $I_{p,c_1}$ is principal then the strict
transform of $\{z_1=0\}$ is
a component of $D_q$. Assume that this is not the case and choose
monomial generators for 
$I_{p,c_1}$ corresponding to lattice points $v_1^*, m_1,\ldots, m_l$ in
$M\cap\sigma^\vee$. We may 
assume that $m_i$ do not contain $v_1^*$, i.e., all $m_i$ lie in the face
$v_1^\perp\cap\sigma^\vee = \langle v_2^*,\ldots, \pm v_n^*\rangle$ of
$\sigma^\vee$. To study the strict transform of $\{z_1=0\}$ in $X_q$ we
first blow up $I_{p,c_1}$, then the rest of the $I_{p,c_i}$, and then
normalize.

Let $Y$ be an affine chart of the blowing up of $X_p$ along $I_{p,c_1}$ (not
necessarily normal), obtained by inverting one of the generators of
$I_{p,c_1}$, and let $E$ be the strict transform of $\{z_1=0\}$ in 
$Y$. Then $E$ is nonempty if and only if  $Y$ is the chart of the
blowing up where we invert one of the $m_i$, say $m_1$. Hence the coordinate
ring 
of $Y$ is generated by monomials corresponding to the lattice points
\[ v_1^*-m_1, m_2-m_1,\ldots,m_l-m_1,v_2^*,\ldots,\pm v_n^*.\]
Since the coefficient of $v_1^*$ in $v_1^*-m_1$ is 1, and the other
generators lie in $v_1^\perp$, we have
\begin{align*}
Y &= \Spec K\left[\frac{z_1}{z^{m_1}},\frac{z^{m_2}}{z^{m_1}},\ldots,
\frac{z^{m_l}}{z^{m_1}}, z_2,\ldots, z_n^{\pm 1}\right] \\
 &= \Spec K\left[\frac{z_1}{z^{m_1}}\right]\times \Spec
K\left[\frac{z^{m_2}}{z^{m_1}},\ldots, 
\frac{z^{m_l}}{z^{m_1}}, z_2,\ldots, z_n^{\pm 1}\right] \\
 &= \bfa^1 \times Y',
\end{align*}
where the strict transform $E$ of $\{z_1=0\}$ is defined by $z_1/z^{m_1}$, on
which $K^*$ acts trivially.

It remains to be shown that if we blow up the ideals $I_{p,c_i}$ for $i\neq 1$
pulled back to $Y$ and normalize, this product structure is preserved. 
We define the ideals $I_{c_i}^Y$ on $Y$ generated by all monomials on which
$K^*$ acts by character $c_i$. The lemma below shows that
$I_{c_i}^Y$ is equal to the inverse image of $I_{p,c_i}$. Hence we may blow
up $I_{c_i}^Y$ instead of the inverse image of $I_{p,c_i}$. Since $K^*$ acts
trivially on $z_1/z^{m_1}$, the ideals $I_{c_i}^Y$ are generated 
by monomials in the second term of the product. Thus, blowing up $I_{c_i}^Y$
preserves the product, and so does normalization. \qed

\begin{lemma}\label{Lem:Torific-preserves-product} For an affine toric variety
$X$ with an action of $K^*$ as a one-parameter subgroup of the torus, let 
$I_\alpha^X$ be the ideal generated by all monomials on which $K^*$ acts by
character $\alpha$. If $\phi: Y\to X$ is a chart of the blowing up of
$I_\alpha^X$ then
\[ I_\beta^Y = (\phi^{-1} I_\beta^X) \cO_Y\]
for all $\beta$.
\end{lemma}

{\bf Proof.} Clearly $\phi^{-1} I_\beta^X\subset I_\beta^Y$.
For the converse, let the monomial generators of the coordinate ring of $Y$
be $z_1/z_{m_1},z^{m_2}/z^{m_1},\ldots,z^{m_l}/z^{m_1},z_1,\ldots,z_n^{\pm 1}$ 
for some generators $z^{m_i}$ of $I_\alpha$. Thus a monomial on $Y$ can be
written as a product 
\[ z^m= (\frac{z_1}{z^{m_1}})^{b_1}(\frac{z^{m_2}}{z^{m_1}})^{b_2}\cdots(\frac{z^{m_l}}{z^{m_1}})^{b_l}
\cdot z_1^{d_1}\cdots z_n^{d_n} \]
for some integers $b_i, d_j\geq 0$ for $i=1,\ldots,l$, $j=1,\ldots,n$. If
$z^m$ happens to be a generator of $I_{\beta}^Y$, i.e., $K^*$ acts on $z^m$
by 
character $\beta$, then also $K^*$ acts on $z^{m'}= z_1^{d_1}\cdots
z_n^{d_n}$ by character $\beta$, and $z^{m'}$ is in
$\phi^{-1}I^X_\beta$. \qed

\begin{corollary}\label{Cor:quotient-toroidal}
The embeddings  $U_{B^{tor}_\pm}/K^* \subset B^{tor}_\pm/K^*$ are toroidal
embeddings, and 
the birational  map $B^{tor}_-/K^* \dra B^{tor}_+/K^*$ is tightly toroidal. 
\end{corollary}

{\bf Proof.} This is immediate from 
the proposition and Lemma~\ref{lem:loc-descr}. \qed

In fact, as
the following lemma, in conjunction with \ref{Lem:Torific-preserves-product},
shows, the map $B^{tor}_-/K^* \dra B^{tor}_+/K^*$ is an isomorphism if the set 
$\{c_1,\ldots,c_\mu\}$ in the definition of the torific ideal $I= I_{c_1}\cdots
I_{c_\mu}$ is chosen large enough. Since we do not need this result, we only
give a sketch of the proof.

\begin{lemma}\label{cor:tor-isom} Let $B=X(N,\sigma) = \Spec
K[z_1,\ldots,z_n^{\pm 1}]$ be a nonsingular affine toric variety, and assume
that $K^*$ acts on $z_i$ by character $c_i$. Let $\alpha\in \ZZ$ be divisible
by all $c_i$, and let $I_\alpha$ and $I_{-\alpha}$ be the ideals generated by
all monomials of $K^*$-character $\alpha$ and $-\alpha$, respectively. If
$\tilde{B}$ is the normalization of the blowing up of $I_\alpha \cdot
I_{-\alpha}$ then the birational map
\[ \tilde{B}_-/K^* \dra \tilde{B}_+/K^* \]
is an isomorphism. The same holds for any torific ideal corresponding to a set
of characters containing $\alpha$ and $-\alpha$.
\end{lemma}

{\bf Sketch of proof.} Let $\sigma = \langle v_1,\ldots,v_m\rangle$, and let
$\pi: 
N_\RR \to N_\RR/\RR\cdot a$ be the projection from $a$. If $\pi$ maps
$\sigma$ isomorphically to $\pi(\sigma)$ then $B_-$ and $B_+$ are isomorphic
already. Otherwise, there exist unique rays $r_+\subset \partial_+ \sigma$
and $r_-\subset \partial_- \sigma$ such that the star subdivision of
$\pi(\partial_+ \sigma)$ at $\pi(r_+)$ is equal to the star subdivision of
$\pi(\partial_- \sigma)$ at $\pi(r_-)$. Now the normalized blowings up of
$I_\alpha$ and $I_{-\alpha}$ turn out to correspond to star subdivisions of
$\sigma$ at 
$r_+$ and $r_-$. The resulting subdivision $\Sigma$ clearly satisfies
$\pi(\partial_- \Sigma) = \pi(\partial_+ \Sigma)$.\qed

It is useful to have a more detailed description of the coordinate ring of
some affine toric charts of $B^{tor}$. The strict transforms of the
divisors $\{z_i = 0\}$ corresponding to the ideals $I_{c_i}$ which are not
principal are removed from the toroidal structure on
$B^{tor}$. Assume $\tau$ is a cone in the subdivision 
associated to the normalization of the blowing up of a torific ideal,
 which contains
$v_1,\ldots,v_k$, the rays in $\tau$ corresponding to the divisors to be
removed from the toroidal structure. We have 
seen above that the 
corresponding affine toric variety $Y$ decomposes as
$$Y = \Spec K[z_i/z^{m_i}] \times Y'.$$
Since $v_j\in \tau$ we have that $(v_i^*-m_i, v_j)\geq 0$ for $i,j=1,\ldots,k$,
Since $m_i$ is positive we have $$(m_i,v_j) = 0, \quad i,j=1,\ldots,k.$$ 
Note that we have a direct product decomposition of cones and lattices dual to
the one in 
condition 4 of Lemma
\ref{Lem:remove-divisors}:
\begin{align*}
\tau^\vee &\cong \langle v_i^*-m_i \rangle \times  \tau' \\
M &\cong M_i \times M_{\widehat\iota},
\end{align*}

 Since,
by condition 3 of 
Lemma \ref{Lem:remove-divisors}, the 
direct product decompositions are compatible,
 we obtain the following:
\begin{corollary}\label{Cor:torific-coordinates}
 Let $B=X(N,\sigma) = \Spec
K[z_1,\ldots,z_n^{\pm 1}]$ be a nonsingular affine toric variety, and assume
that $K^*$ acts on $z_i$ by character $c_i$.
Let $Y \subset B^{tor}$ be an
affine toric chart corresponding to a cone $\tau$ containing
$v_1,\ldots,v_k$, the rays in $\tau$ corresponding to the divisors to be
removed from the toroidal structure. Then there exist $m_i \in \sigma^\vee$
such that $(m_i,v_j) = 
0 $ for $i,j=1,\ldots,k$
 and $z_i/z^{m_i}$ are
invariant, and a toric variety $Y'$,
 such that 
$$Y = \Spec K\left[\frac{z_1}{z^{m_1}},\ldots, 
\frac{z_k}{z^{m_k}}\right] \times Y'.$$
\end{corollary}

\begin{example}
Consider $B= {\bbA}^3 = \Spec K[z_1,z_2,z_3]$, where $t \in K^*$ acts as
$$t\cdot(z_1, z_2, z_3) = (t^2z_1, t^3z_2, t^{-1}z_3).$$ 
We have the following generators of the torific ideals
$I_\alpha$:
\begin{align*}
I_2 &= \{z_1, z_2z_3\} \\
I_3 &= \{z_2, z_1^2z_3\} \\
I_6 &= \{z_1^3, z_2^2, z_1^2z_2z_3\}\\
I_{-1} &= \{z_3\}. 
\end{align*}
Let $I = I_2 I_3 I_6 I_{-1}$. 
If we regard $B = X(N,\sigma)$ as the toric variety corresponding
to the  cone 
$$\sigma = \langle v_1,v_2,v_3\rangle \subset N_{\mathbb R},$$
then $B^{tor}$
is described by the fan covered by the following four maximal cones
$$\begin{array}{rl}
\sigma_1 &= \langle v_1,v_1+v_3,v_1 + v_2\rangle \\
\sigma_2 &= \langle v_1 + v_2,v_1+v_3, 2 v_1 + 3 v_2,v_3\rangle\\
\sigma_3 &= \langle 2 v_1+3 v_2, v_3,2v_2 + v_1,v_1+v_2,v_2 + v_3\rangle.\\
\sigma_4 &= \langle 2 v_2+v_1, v_2+v_3, v_2 \rangle 
\end{array}
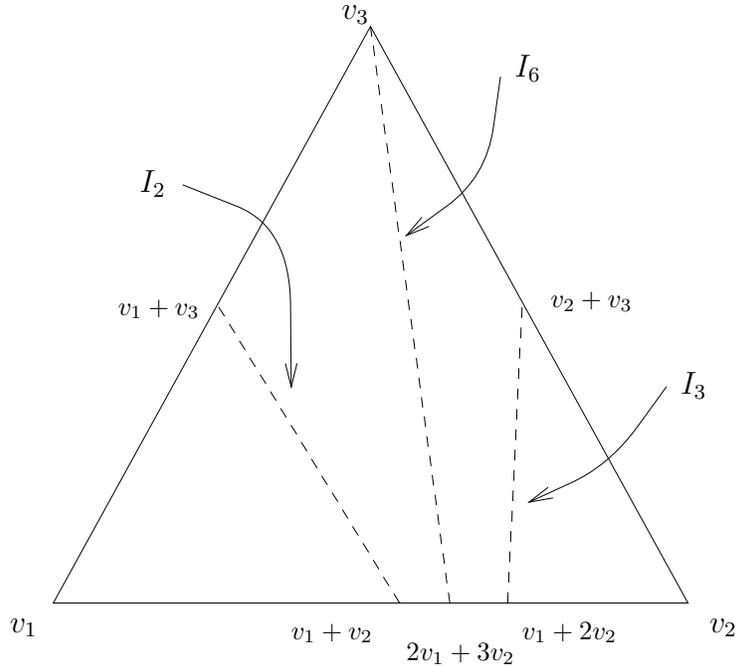
$$
\begin{figure}[htb]
\begin{center}
\setlength{\unitlength}{0.0005in}
\begingroup\makeatletter\ifx\SetFigFont\undefined
% extract first six characters in \fmtname
\def\x#1#2#3#4#5#6#7\relax{\def\x{#1#2#3#4#5#6}}%
\expandafter\x\fmtname xxxxxx\relax \def\y{splain}%
\ifx\x\y   % LaTeX or SliTeX?
\gdef\SetFigFont#1#2#3{%
  \ifnum #1<17\tiny\else \ifnum #1<20\small\else
  \ifnum #1<24\normalsize\else \ifnum #1<29\large\else
  \ifnum #1<34\Large\else \ifnum #1<41\LARGE\else
     \huge\fi\fi\fi\fi\fi\fi
  \csname #3\endcsname}%
\else
\gdef\SetFigFont#1#2#3{\begingroup
  \count@#1\relax \ifnum 25<\count@\count@25\fi
  \def\x{\endgroup\@setsize\SetFigFont{#2pt}}%
  \expandafter\x
    \csname \romannumeral\the\count@ pt\expandafter\endcsname
    \csname @\romannumeral\the\count@ pt\endcsname
  \csname #3\endcsname}%
\fi
\fi\endgroup
{\renewcommand{\dashlinestretch}{30}
\begin{picture}(7721,6918)(0,-10)
\path(450,636)(3750,6636)(7050,636)
	(450,636)(450,636)
\dashline{120.000}(2175,3711)(4050,636)
\dashline{120.000}(5325,3711)(5175,636)
\dashline{120.000}(3750,6636)(4575,636)
\spline(1800,4986)
(2925,4536)(2925,2886)
\path(2865.000,3126.000)(2925.000,2886.000)(2985.000,3126.000)
\spline(5100,6111)
(4950,5061)(4125,4461)
\path(4283.806,4650.685)(4125.000,4461.000)(4354.387,4553.637)
\spline(6825,2886)
(6225,2061)(5400,1686)
\path(5593.660,1839.935)(5400.000,1686.000)(5643.316,1730.691)
\put(0,336){\makebox(0,0)[lb]{\smash{{{\SetFigFont{12}{14.4}{rm}$v_1$}}}}}
\put(3450,6711){\makebox(0,0)[lb]{\smash{{{\SetFigFont{12}{14.4}{rm}$v_3$}}}}}
\put(7275,336){\makebox(0,0)[lb]{\smash{{{\SetFigFont{12}{14.4}{rm}$v_2$}}}}}
\put(1125,3636){\makebox(0,0)[lb]{\smash{{{\SetFigFont{10}{12}{rm}$v_1+v_3$}}}}}
\put(5625,3711){\makebox(0,0)[lb]{\smash{{{\SetFigFont{10}{12}{rm}$v_2+v_3$}}}}}
\put(2925,261){\makebox(0,0)[lb]{\smash{{{\SetFigFont{10}{12}{rm}$v_1+v_2$}}}}}
\put(4125,36){\makebox(0,0)[lb]{\smash{{{\SetFigFont{10}{12}{rm}$2v_1+3v_2$}}}}}
\put(5325,261){\makebox(0,0)[lb]{\smash{{{\SetFigFont{10}{12}{rm}$v_1+2v_2$}}}}}
\put(1350,4911){\makebox(0,0)[lb]{\smash{{{\SetFigFont{12}{14.4}{rm}$I_2$}}}}}
\put(5250,6111){\makebox(0,0)[lb]{\smash{{{\SetFigFont{12}{14.4}{rm}$I_6$}}}}}
\put(6975,2811){\makebox(0,0)[lb]{\smash{{{\SetFigFont{12}{14.4}{rm}$I_3$}}}}}
%\put(5400,5961){\makebox(0,0)[lb]{\smash{{{\SetFigFont{10}{12}{rm}$6$}}}}}
%\put(7125,2661){\makebox(0,0)[lb]{\smash{{{\SetFigFont{10}{12}{rm}$3$}}}}}
%\put(1500,4761){\makebox(0,0)[lb]{\smash{{{\SetFigFont{10}{12}{rm}$2$}}}}}
\end{picture}
}
%\centerline{\psfig{figure=exa.ps,width=3in}}
\end{center}
\end{figure}

The dual cone $\sigma_1^{\vee}$ has the product description
$$\begin{array}{rl}
\sigma_1^{\vee} &= \langle v_1^* - (v_2^* + v_3^*), v_2^*, v_3^*\rangle \\
&= \langle v_1^* - (v_2^* + v_3^*)\rangle \times \langle v_2^*,
v_3^*\rangle. 
\end{array}$$
Thus, even if
we remove the divisor
$\{z_1/z_2z_3 = 0\}$ from  the original toric
structure of $$X(N,\sigma_1) = \Spec k[z_1/z_2z_3,z_2,z_3],$$ we still have the
toroidal embedding structure 
$$X(N,\sigma_1) \setmin (\{z_2 = 0\} \cup \{z_3 = 0\}) \subset X(N,
\sigma_1). $$ 
As $z_1/z_2z_3$ is invariant, the action of $K^*$ is toroidal.
For example, at $0\in X(N, \sigma_1)$ we have a toric chart
\begin{align*}
K^* \times K^2 &\to K\times K^2\cong X(N,\sigma_1)\\
(x_1,x_2,x_3) &\mapsto (x_1-1,x_2,x_3).
\end{align*}
Globally, the divisors corresponding to the new
rays 
$$D_{\langle v_1 + v_2\rangle},D_{\langle v_1 + v_3\rangle},D_{\langle 2v_1 +
3 v_2\rangle}, D_{\langle v_1 + 2 v_2\rangle}, D_{\langle v_2 + v_3\rangle}$$
together with $D_{\langle v_3\rangle}$ coming from $I_{-1}$, are obtained
through the blowing up of the torific ideals. Considering 
$$ U_{B^{tor}} = B^{tor} \setmin (D_{\langle v_1 + v_2\rangle}\cup D_{\langle
v_1+ v_3\rangle}\cup D_{\langle 2v_1 + 3 v_2\rangle}\cup D_{\langle v_1 + 2
v_2\rangle} \cup D_{\langle v_2 + v_3\rangle} \cup D_{\langle v_3\rangle})$$
we obtain a toroidal structure $U_{B^{tor}} \subset B^{tor}$  with
a toroidal $K^*$-action. 
\end{example}

\section{A proof of the weak factorization theorem}\label{Sec:connecting}

\subsection{The situation}
In Theorem \ref{Th:locally-toric-factorization}  we have
constructed a factorization of 
the given birational map $\phi$ into tightly locally toric birational maps
\[
\begin{array}{rclcrclcrcl}
X_1 = W_{1-}  & \dra &  W_{1+}& \cong & W_{2-} &\dra & W_{2+}& \ldots & W_{m-}
&\dra & 
W_{m+} = X_2, \\ 
\searrow & &\swarrow & &  \searrow & &\swarrow & & \searrow & &\swarrow \\
& B_{a_1}/\!/ K^* & & & & B_{a_2}/\!/ K^* & & & & B_{a_m}/\!/ K^*& 
\end{array}
\] 
where $W_{i\pm} = (B_{a_i})_\pm/K^*$ (here $W_{i-}$ is $W_{i-1}$ in the
notation of 
Theorem 
\ref{Th:locally-toric-factorization}, and $W_{i+}$ is $W_{i}$).

For a choice of a torific ideal $I= I_{c_1}\cdots
I_{c_\mu}$ on $B_{a_i}$, denote by $B_{a_i}^{tor}
\rightarrow B_{a_i}$  the corresponding torific blowing up.  Write 
$W_{i\pm}^{tor} = B_{a_i\pm}^{tor}/K^*$, and $U_{i\pm}^{tor} =
U_{B_{a_i\pm}^{tor}}/K^*$. We have a natural diagram of birational maps 
\[
\begin{array}{ccccc}
W_{i-}^{tor} & & \dra & & W_{i+}^{tor} \\
\dar f_{i-}& \searrow & & \swarrow & \dar f_{i+} \\
W_{i-} & & B_{a_i}^{tor}/\!/ K^* & & W_{i+} \\
 &  \searrow & \dar & \swarrow & \\
 & & B_{a_i}/\!/ K^* & &
\end{array}
\]
By Corollary \ref{Cor:quotient-toroidal} the  embeddings  $U_{i\pm}^{tor}
\subset  W_{i\pm}^{tor}$ are toroidal, and the birational map  
$\varphi_i^{tor}: W_{i-}^{tor} \dra W_{i+}^{tor}$ is tightly toroidal.

We say that the ideal $I = I_{c_1}\cdots
I_{c_\mu}$ is {\em balanced} if $\sum c_j = 0$. It follows from lemma
\ref{Lem:ideal-nonzero} that we can always enlarge the set
$\{c_1,\ldots,c_\mu\}$ 
to get a  balanced torific ideal $I$.

\begin{lemma} Suppose the torific ideal $I$ is balanced. Then the morphism
$f_{i\pm}$ is a blowing up of a canonical ideal sheaf 
$I_{i\pm}$ on $W_{i\pm}$.  
\end{lemma}

{\bf Proof.} By 
Lemma~\ref{Lem:ideal-existance-toric}, the ideal $I$ is generated by
$K^*$-invariant sections, and we can identify $I$ as the inverse image of 
an ideal sheaf in $B_{a_i}/\!/K^*$ generated by the same sections. Let
$I_{i\pm}$ be the pullback of this ideal sheaf to $(B_{a_i})_\pm/K^*$ via the 
map $(B_{a_i})_\pm/K^* \to B_{a_i}/\!/K^*$. Then $f_{i\pm}$ is the blowing up
of 
$I_{i\pm}$ because taking the quotient by $K^*$ commutes with blowing up the
sheaf $I$. \qed

We note that this lemma is true even when $I$ is not balanced; however if $I$
is not balanced the construction of a {\em canonical} ideal sheaf is less
immediate. 
{\em From now on we assume that the torific ideals are chosen to be balanced.}

Note that if the varieties $W_{i\pm}$ were nonsingular and the morphisms
$f_{i\pm}$ were composites of blowings up of smooth centers, we would get the
weak factorization by applying Theorem~\ref{thm:morelli} to each 
$\varphi_i^{tor}$. This is not the case in general.  In this section we replace
$W_{i\pm}$ by nonsingular varieties and 
$f_{i\pm}$ by composites of blowings up with nonsingular centers.

\subsection{Lifting toroidal structures}\label{Sec:lifting-toroidal}
Let $W_{i\pm}^{res} \rightarrow W_{i\pm}$ be the canonical resolution of
singularities.  Note that, since $W_{i+} = W_{(i+1) -}$, we have 
$W_{i+}^{res} = W_{(i+1) -}^{res}$. 

Denote $I_{i\pm}^{res} = I_{i\pm} \cO_{W_{i\pm}^{res}}$. Let $W_{i\pm}^{can}
\to W_{i\pm}^{res}$ be the canonical principalization of the ideal
$I_{i\pm}^{res}$, and let $h_{i\pm}: W_{i\pm}^{can} \to W_{i\pm}^{tor}$ be the
induced morphism. 

\[ 
\begin{array}{cclcccrcc}
W_{i-}^{can} & \stackrel{h_{i-}}{\to} & W_{i-}^{tor} &  & \dra  &  &
W_{i+}^{tor} & \stackrel{h_{i+}}{\leftarrow} & W_{i+}^{can}\\ 
\dar & & \dar f_{i-} & \searrow & & \swarrow &  \dar f_{i+} & & \dar \\
W_{i-}^{res} & \to & W_{i-} &  & B_{a_i}^{tor}/\!/ K^*  &  & W_{i+} &
\leftarrow & W_{i+}^{res}\\ 
& & & \searrow & \dar & \swarrow & & & \\
& & & & B_{a_i}/\!/ K^*& & & &
\end{array}
\]

 Denote $U_{i\pm}^{can} =h_{i\pm}^{-1}  U_{i\pm}^{tor}$.
The crucial point now is to show:

\begin{proposition} The embedding $U_{i\pm}^{can} \subset W_{i\pm}^{can}$ 
is a toroidal embedding, and the morphism $W_{i\pm}^{can} \to
W_{i\pm}^{tor}$ is toroidal.
\end{proposition}

{\bf Proof.} For simplicity of notation we drop the subscripts $i$ and $a_i$,
as we treat 
each quasi-elementary piece separately. We may assume that all the varieties
$B, W_{\pm}, 
W_{\pm}^{tor}, W_{\pm}^{res}, W_{\pm}^{can}$ and the morphisms between
them are toric. Indeed, if $V_p\to X_p$ is a toric chart at some point
$p\in W_{\pm}$, obtained from a toric chart in $B$, we get a
toric chart for $W_{\pm}^{tor}$ by blowing up a 
torific ideal in $X_p$, which is a toric ideal since it is generated by
monomials. Similarly, resolution of singularities and 
principalization over the toric variety $X_p$ provide toric charts for 
$W_{\pm}^{res}$ and $W_{\pm}^{can}$. The maps are toric (i.e., torus
equivariant) by canonicity. 

Consider now the diagram of toric morphisms between toric varieties and the
corresponding diagram of fans:
\[ \begin{array}{ccccccc}
W_{\pm}^{can} & \to & W_{\pm}^{tor} & \hspace{.5in} & \Sigma_{\pm}^{can} & 
\to & \Sigma_{\pm}^{tor} \\
& \searrow & \dar & & & \searrow & \dar \\
& & W_{\pm} & & & & \Sigma_{\pm}
\end{array}
\]
Let $X_\sigma\subset W_{\pm}^{tor}$ be an affine open toric subvariety
corresponding to a cone $\sigma\in \Sigma_{\pm}^{tor}$, and 
write
\[ X_\sigma \cong \bfa^k \times X_{\sigma'} ,\]
where the toric divisors $E_1,\ldots,E_k$ pulled back from $\bfa^k$ are the
ones 
removed in order to define the 
toroidal structure on $W_{\pm}^{tor}$. Let $X_\sigma^{can}$ be the inverse
image of $X_\sigma$ in $W_{\pm}^{can}$. We need to show that we have a
decomposition $X_\sigma^{can}\cong \bfa^k \times X^{can}_{\sigma'},$  such that
the resulting map  $\bfa^k \times X^{can}_{\sigma'} \to \bfa^k \times
X_{\sigma'}$  is a product, with the second factor being the identity map.

Write $X_\sigma = B_\sigma/K^*$, where $B_\sigma \subset
B^{tor}_\pm$ is the affine open toric subvariety lying over 
$X_\sigma$.   

 By Corollary \ref{Cor:torific-coordinates}, the coordinate  
rings of $B_\sigma$ and $X_\sigma$ can be written as
\begin{align*}
 A_{X_\sigma} &\cong
K[\frac{z_1}{z^{m_1}},\ldots,\frac{z_k}{z^{m_k}}]\otimes  A_{X_{\sigma'}},\\
 A_{B_\sigma} &\cong
K[\frac{z_1}{z^{m_1}},\ldots,\frac{z_k}{z^{m_k}}]\otimes  A_{B_{\sigma'}},
\end{align*}
where $X_{\sigma'} = B_{\sigma'}/K^*$, and
where $z^{m_j}$ are monomials on which $K^*$ acts with the same character as
on $z_j$, such that $z_i\nmid z^{m_j}$ for $i,j=1,\ldots,k$.

\begin{lemma} For each $y=(y_1,\ldots,y_k) \in K^k$ consider the automorphism
$\theta_y$ of  $B$ defined by 
$$\begin{array}{rll}
\theta_y(z_i) &= z_i + y_i \cdot  z^{m_i}, & i\leq k\\
 \theta_y(z_i) &= z_i, & i>k.
\end{array}$$
Then
\begin{enumerate}
\item $\theta_y$ defines an action of the additive group $K^k$ on $B$.
\item The action of $\theta_y$ commutes with the given $K^*$-action.
\item The ideals  $I_{c}$ are invariant under this action.
\item The action leaves $B_\pm$ invariant, and descends to $W_\pm$.
\item The action lifts to $B^{tor}$.
\item This action on $B^{tor}$  leaves the open set  $B_\sigma$ invariant.
\item The induced  action on  $B_\sigma$ descends to a fixed-point-free action
of $K^k$ on 
$X_\sigma$. 
\item The resulting action on $X_\sigma$ is given by
$$\bar\theta_y(z_i/m_i) = z_i/m_i+ y_i; \quad \theta_c(f) = f \text{\ for\
} f\in  A_{X_{\sigma'}}.$$
\end{enumerate}
\end{lemma}
{\bf Proof.} Since $z_i\nmid m_j$  for $i,j=1,\ldots,k$, we have that
the $\theta_y$ commute with each other, and  $\theta_y + \theta_{y'}
=\theta_{y+y'}$ thus defining a $K^k$-action. Since 
$K^*$ acts on $z_i$ and $m_i$ through the same character, it 
commutes with $\theta_y$. For the same reason the ideals $I_c$ are
invariant. Since $B_- = B \setmin V(\sum_{c<0} I_{c})$ we have that $B_-$ is
invariant, and similarly for $B_+$; since the $K^k$-action commutes with $K^*$
it  
descends to $W_\pm$. Since $I = \prod I_{c_i}$ we have that $I$
is $K^k$-invariant and therefore the $K^k$-action lifts to  $B^{tor}$. Also 
by definition $\theta_y(z_i/m_i) = z_i/m_i+ y_i$, which implies the rest of the
statement. \qed

Back to the proposition. Since $W_\pm^{res} \to W_\pm$ is the canonical
resolution of singularities, the action of $K^k$ lifts to $W_\pm^{res}$. Since
the ideal $I_\pm$ is generated by $K^*$-invariants in $I$, and since the action
of $K^*$
commutes with $\theta_c$, we have that $I_\pm$ is invariant under $K^k$, and
therefore $I_\pm^{res}$ is invariant under $K^k$ as well. Since $W_\pm^{can}
\to W_\pm^{res}$ is the canonical principalization of $I_\pm^{res}$, the action
of  $K^k$ lifts to $W_\pm^{can}$. In particular, the map  $W_\pm^{can} \to
W_\pm^{tor}$ is $K^k$-equivariant. By the lemma, the action of $K^k$ on the
invariant open set $X_\sigma\subset W_\pm^{tor}$ is fixed-point free, therefore
the action on the inverse image   $X_\sigma^{can}$ is  fixed-point
free. Writing $X_{\sigma'}^{can}$ for the inverse image of 
$(0,\ldots,0)\times X_{\sigma'}$, we have an equivariant  decomposition
$W_\pm^{can}\cong \bfa^k \times
X_{\sigma'}^{can}$ as needed. \qed

\subsection{Conclusion  of the proof}

Since $X_1 = W_{1-}$ and $X_2 = W_{m+}$ are nonsingular, we have $W_{1-}^{res}
= 
W_{1-}$ and $W_{m+}^{res} = W_{m+}$.  For each $i = 1, \ldots, m$ we
have obtained a diagram

$$\begin{CD}
@.@.@.W_{i-}^{can} @.
\overset{\varphi_i^{can}}{\dashrightarrow} @.@.W_{i+}^{can} @.@.@. \\
@. r_{i-} \swarrow \hskip.2in @.@.@VV{h_{i-}}V @.@.@VV{h_{i+}}V
@.@. \hskip.2in \searrow r_{i+}
\\
W_{i-}^{res}@.@.@.W_{i-}^{tor} @.
\overset{\varphi_i^{tor}}{\dashrightarrow} @.@.W_{i+}^{tor}
@.@.@.@. W_{i+}^{res}
\\
@VVV \hskip.1in {f_{i-}} \swarrow \hskip.2in @.@.@.@.@.@.@.@. \hskip.2in
\searrow {f_{i+}} \hskip.1in @.@VVV
\\
W_{i-} @.@.@.@.\overset{\varphi_i}{\dashrightarrow} @.@.@.@.@.@. W_{i+}
\\
\end{CD}$$
where

\begin{enumerate}\item
 the canonical principalizations $r_{i-}$ and $r_{i+}$
are composites of blowings up with smooth centers,
\item
  $\varphi_i^{can}$ is tightly toroidal.
\end{enumerate}

Applying Theorem~\ref{thm:morelli-proj} to the toroidal map $\varphi_i^{can}$ 
we see that $\varphi_i^{can}$ is a composite of 
toroidal blowings up and blowings down,
with smooth centers, between nonsingular toroidal embeddings. Thus we get a 
factorization 
$$\phi:X_1 = W_{1-}^{res} \dra W_{1+}^{res}=W_{2-}^{res} {\dashrightarrow}
\cdots \dashrightarrow
W_{m-}^{res} {\dashrightarrow} W_{m+}^{res} = X_2,$$
where all $W_i^{res}$ are nonsingular, and the birational maps are
composed of a sequence of blowings up and blowings
down. We do not touch the open subset $U\subset X_1$ on which $\phi$ is an
isomorphism. After the reduction step in Lemma \ref{lem:red-proj}, the 
projectivity over $X_2$ follows from Proposition~\ref{prop:projectivity}, the
projectivity statement in Theorem~\ref{thm:morelli-proj}, and
the construction. Finally, blowing up a nonsingular center can be
factored as a sequence of blowings up of {\em irreducible} centers, simply
blowing up 
one  connected component  at a time; since blowing up is a projective
operation, this preserves projectivity.  This completes the proof of 
Theorem~\ref{Th:weak-factorization}. \qed

\section{Generalizations}\label{Sec:generalizations}

\subsection{Reduction to an algebraically closed overfield}
We begin our proof of Theorem \ref{Th:general-weak-factorization}. We claim
that, 
in case (1) of algebraic spaces, it suffices to prove the result in case $L$ is
algebraically closed. Let $\bar L$ be an algebraically closed field containing
$L$. Given $\phi:X_1 \das X_2$, isomorphic on $U$,  consider
the map $\phi_{\bar L}:(X_1)_{\bar L} \das 
(X_2)_{\bar L}$.  Assuming the generalized factorization theorem applies over
such a field, we get ${\varphi_i}_{\bar L}: \bar V_i \das \bar V_{i+1}$. The
functoriality of this  
factorization guarantees that the Galois group acts on $\bar V_i$, and
${\varphi_i}_{\bar L}$ are Galois equivariant. Therefore, denoting $V_i = \bar
V_i/Gal(\bar L/L)$, we get $\varphi_i:V_i\das  V_{i+1}$ as required.

\subsection{Reduction to an algebraically closed subfield} Still considering
case (1),  suppose $L\subset K$ are algebraically closed fields, and suppose we
have the theorem for algebraic spaces over fields isomorphic to $L$. If
$\phi:X_1 \das X_2$ is a birational map over $L$, with factorization given by
$\varphi_i:V_i\das  V_{i+1}$, then we claim that the induced maps
${\varphi_i}_K: {V_i}_K \das  {V_{i+1}}_K$ is functorial over $K$. Indeed,
any isomorphism $K \to K'$ carries $L$ to an isomorphic field, and the
functoriality over $L$ induces the desired morphisms ${V_i}_K \to
{V_i}'_{K'}$.

\subsection{Reduction to $L=\CC$} Still considering
case (1), let $K$ be algebraically closed and let $\phi:X_1 \das X_2$ be a
birational map of complete algebraic spaces over $K$. Then, by definition,
$X_i$ are given by  \'etale equivalence relations $R_i \subset Y_i^2$, where
$R_i$ and $Y_i$ are varieties over $K$, and $\phi$ is defined by suitable
correspondences between $Y_i$. Also the open set $U$ corresponds to a Zariski
open in $Y_i$. All these varieties can be defined over a 
finitely generated subfield $L_0 \subset K$, and therefore over its algebraic
closure $L \subset K$. But any such $L$ can be embedded in $\CC$. Therefore, by
the previous reductions, it suffices to consider the case of algebraic spaces
over a field $L$ isomorphic to $\CC$.

By considering the associated analytic spaces, this allows us to use structures
defined in the analytic category, as long as we ensure that the resulting
blowings up are functorial in the algebraic sense, namely, independent of a
choice of isomorphism $L \to \CC$.

\subsection{Reduction to a projective morphism}
Now we consider both cases (1) and (2). To simplify the terminology, we use the
term ``birational map'' to indicate also a bimeromorphic map.
Given $\phi:X_1 \das X_2$ isomorphic on $U$, let $X'_i \to X_i$ be the
canonical principalizations of $X_i\setmin U$ (endowed
with reduced structure).  It is convenient to replace $X_i$ by $X_i'$ and
assume from now on that $X_i\setmin U$ is a simple normal crossings divisor.

 We note
 that   Lemma  \ref{lem:red-proj} works word for word in the cases of
 algebraic spaces or analytic spaces. As we have already remarked, this
 procedure is functorial. Also, the centers
 of blowing  up have normal crossings with the inverse image of $X_i \setmin
 U$. 

It is also easy to see that the resulting morphism $X'_1 \to X'_2$ is endowed
with a relatively ample line bundle which is functorial under absolute
isomorphisms. Indeed, the $\Proj$ construction of a blowing up gives a
functorial relatively ample line bundle for each blowing up. Furthermore, if
$f_1:Y_1  \to Y_2 $ and  $Y_2  \to Y_3 $ are given relatively ample line
bundles  $L_1$ and $L_2$, then there is a minimal positive integer $k$ such
that $L_1 
\otimes 
f_1^*L_2^{\otimes k}$ is relatively ample for $Y_1 \to Y_3$; thus we can form a
functorial 
relatively ample 
line bundle for a sequence of blowings up. In an analogous manner we can form
a functorial  ideal sheaf $I$ on  $X'_2$ such that $X_1'$ is the blowing up of
$I$. 

From now on we assume $X_i\setmin U$ is a simple normal crossings divisor and
$\phi$ is a projective morphism.

\subsection{Analytic locally toric structures} 
There are various settings in which one can generalize locally toric and
toroidal structures to algebraic and analytic spaces, either using formal
completions (see \cite{KKMS}), or \'etale charts (see \cite{Matsuki-notes}), or
logarithmic structures (see \cite{Kato}). Here 
we try to keep things simple, by sticking to the analytic situation, and
modifying our earlier definitions slightly.

An {\em analytic toric chart} $V_p \subset W$, $\eta_p: V_p \to X_p$  is
defined to be a neighborhood  of $p$ in the euclidean topology, with $\eta_p$
an {\em open immersion} in the euclidean topology. The fact that we use open
immersions simplifies our work significantly.

The notions of analytic locally toric structures, analytic toroidal embeddings,
modifications, toroidal birational maps and tightly toroidal  birational maps
are defined as in the  
case of varieties, using analytic toric charts. 

We note that in an analytic toroidal embedding, the toroidal divisors may have
self intersections. If $U\subset X$  is an analytic toroidal embedding, and if
$X'\to X$ is the canonical embedded resolution of singularities of $X \setmin
U$, then $X'\setmin U$ is a {\em strict toroidal embedding}, namely one without
self intersections.

For strict toroidal embeddings, the arguments of \cite{KKMS} regarding rational
conical 
polyhedral complexes, modifications and subdivisions go through, essentially
word for word. The divisorial description of the cones (see
\cite{KKMS}, page 61) shows that the association $(U\subset X) \mapsto
\Delta_X$ of a polyhedral complex to a toroidal embedding is
functorial under absolute isomorphisms in both the analytic and algebraic
sense, and similarly for the modification  associated to a subdivision.  

\subsection{Functorial toroidal factorization} 
Consider an analytic toroidal birational map $\phi: W_1 \das W_2$ of complete
nonsingular toroidal embeddings $U\subset W_i$. By the
resolution of singularities argument above, we may assume $U \subset W_i$ are
strict toroidal embeddings. Theorem \ref{thm:morelli-proj} applies in
this situation, but we 
need to make the construction functorial. It may be appropriate to rewrite
the proof in a functorial manner, but this would take us beyond the
intended scope of this paper. Instead we show here that the result can
be made equivariant under the automorphism group of a fan cobordism, which,
assuming the axiom of choice, implies functoriality.

Let $\Delta_i$ be the polyhedral complex of $U \subset W_i$.  
  Denote by $G_i$ the automorphism group of $\Delta_i$. Since
an automorphism of   $\Delta_i$ is determined by its action on the primitive
points of the rays in $\Delta_i$, these groups are finite. 

Consider the barycentric subdivision $B\Delta_i\to\Delta_i$ (see \cite{KKMS},
III 2.1, or \cite{Abramovich-Wang}). It corresponds  to a composition of 
blowings up $BW_i \to W_i$, which is functorial. The group $G_i$ acts on
$B\Delta_i$. The subdivision   $B\Delta_i\to\Delta_i$ has the following
property: given a cone $\sigma$ in $B\Delta_i$, an element $g\in G_i$, and a
ray $\tau$ in $\sigma$ such that $g\tau$ is also in $\sigma$, we have
$g\tau=\tau$. This means, in particular, that for any subgroup $H\subset G_i$
and any $H$-equivariant subdivision $\Delta\to B\Delta_i$
the quotient $\Delta/H$ is also a polyhedral complex (see
\cite{Abramovich-Wang}).

Let $Z$ be the canonical resolution of singularities of the graph of $BW_1
\das BW_2$. This is clearly functorial in $\phi$. Now $Z \to BW_i$ are
toroidal birational morphisms, corresponding 
to  subdivisions $\Delta_Z \to B\Delta_i$. Let $H\subset G_1$ be the subgroup
stabilizing the subdivision $\Delta_Z \to B\Delta_1$. 

Fix a representative in the isomorphism class of $\Delta_Z \to B\Delta_1$, and,
using the axiom of choice, fix an isomorphism of any element of the isomorphism
class with this representative.  
Note that the absolute automorphism
group of $Z\to W_1$ maps to $H$. Therefore, in order to construct a functorial
factorization 
of $Z\to W_1$ it suffices to construct an $H$-equivariant combinatorial
factorization of our representative of the isomorphism class, which by abuse of
notation we call $\Delta_Z \to
B\Delta_1$.

Now $\Delta_Z/H \to   
B\Delta_1/H$ is a subdivision of nonsingular polyhedral complexes, and
the toroidal weak factorization theorem says that it  
admits a combinatorial factorization, as a sequence composed of nonsingular
star 
subdivisions and inverse nonsingular 
star subdivisions. Lifting these subdivisions to   $\Delta_Z
\to B\Delta_1$,  
we get the resulting $H$-equivariant factorization, which in turn corresponds
to a toroidal factorization of $BW_1 \das Z$. We now apply the
same procedure to $Z' \to BW_2$. This gives the desired  functorial toroidal
factorization of $\phi$. 

\subsection{Analytic toroidal $\CC^*$-actions} The nature of  $\CC^*$-actions
on analytic spaces differ significantly from the case of varieties. 
However, the situation is almost the same if one restricts to {\em relatively
algebraic} actions.

\begin{definition} Let $X \to S$ be a morphism of analytic spaces and $L$ a
relatively ample 
line bundle for $X\to S$. An action of
$\CC^*$ on $X,L$ over $S$ is {\em relatively algebraic} if there is an open
covering $S = \cup  
S_i$, an algebraic action of $\CC^*$ on a projective space $\PP^{N_i}$, and a
Zariski-locally-closed $\CC^*$-equivariant embedding $X \times_SS_i \subset
S_i\times  \PP^{N_i}$, such that for some integer $l_i$ we have that
$L^{l_i}_{X \times_SS_i}$ is $\CC^*$-isomorphic to the pullback of
$\cO_{\PP^{N_i}}(1)$. 
\end{definition} 

It is easy to see that if $X \to S$ is a projective morphism, $L$ a line
bundle, 
with a relatively 
algebraic $\CC^*$-action, then $X \subset {\cP}roj_S \Sym E$, {\em where the
sheaf $E
=\oplus_{i=1}^k E_i$ is a completely reducible $\CC^*$ sheaf.}

In the analytic category we use embedded charts rather than \'etale
ones. Accordingly, we say that a  $\CC^*$-equivariant open set $V \subset X$ is
{\em strongly embedded} if for any orbit $O\subset V$, the closure of $O$ in
$X$ is contained in $V$. This implies that $V/\!/\CC^* \to X/\!/\CC^*$ is an
open embedding. We define an analytic locally toric $\CC^*$-action on $W$
using strongly embedded toric charts $\eta_p: V_p \to X_p$ (we still have the
requirement that $V_p= \pi^{-1} \pi V_p$, where $\pi:W\to W/\!/K^*$ is the
projection, which means that $V_p \subset W$ is also strongly embedded).  

It is not difficult to show that a strongly embedded toric chart exists for
each point $p\in B$, the analogue of Luna's fundamental lemma.

With these modification, Lemma \ref{lem:loc-tor-act} is proven in the same
manner in the analytic setting. We also note that, if $D=\sum_{i=1}^l D_i
\subset W$ is 
a simple normal crossings divisor, then toric charts can be chosen compatible
with $D$. Indeed, we only need to choose semi-invariant parameters
$x_1,\ldots,x_n$ so that $x_i$ is a defining equation for $D_i$, for
$i=1,\ldots,l$. 

\subsection{Analytic birational cobordisms} Analytic birational cobordisms are
defined the same way as in the case of varieties, with the extra assumption
that the $\CC^*$-action is relatively algebraic. 

Given a projective birational
morphism $\phi:X_1 \to X_2$ we construct a compactified, relatively projective
cobordism $\overline B \to X_2$ as in the algebraic situation, with the
following modification: using canonical resolution of singularities we make the
inverse image of $X_2\setmin U$ in 
$\overline B$ into a simple normal crossings divisor, crossing $X_1$ and $X_2$
normally. Note that these operations are functorial in absolute isomorphisms of
$\phi$. 

As indicated before, this construction endows $\overline B\to X_2$ with a
functorial relatively ample line bundle. Since this bundle is obtained from the
$\Proj$ construction of the blowing up of an invariant ideal, it comes with a
functorial $\CC^*$-action as well.

The
considerations of collapsibility and geometric invariant theory work as in the
algebraic setting, leading to  Theorem
\ref{Th:locally-toric-factorization}. We note that the resulting locally toric
factorization is functorial, and  the toric charts on $W_i$ can be chosen
compatible with the divisor coming from $X_1\setmin U$ or $X_2\setmin U$.

\subsection{Functoriality of torification and compatibility with divisors}
We note that the definition of the torific ideals is clearly functorial. The
proof of its existence works as in the case of varieties. The same is true for
its torifying property. In order to make this constrcuction compatible with
divisors, we replace the total transform $D$ of $I$ by adding the inverse image
of $X_2\setmin U$. This guarantees that the resulting toroidal structure on
$B^{tor}$ is compatible with the divisors coming from  $X_2\setmin U$.

\subsection{Conclusion of proof of Theorem \ref{Th:general-weak-factorization}}
Canonical resolution of singularities is functorial, therefore the construction
of $W^{res}_\pm \to W_\pm$ is functorial. We can now replace $W^{res}_\pm$ by
the 
canonical principalization of the inverse image of $X_2\setmin U$, making the
latter a {\em simple} normal crossings divisor. 
 Since the ideal $I$ is functorial, the construction of $W^{can}_\pm \to W_\pm$
is functorial, and the locally toric structure implies that the centers of
blowing up in  $W^{can}_\pm \to W^{res}_\pm$ have normal crossings with  the
inverse image of $X_2\setmin U$. We can now apply functorial toroidal
factorization  to the toroidal birational map $W^{can}_- \das
W^{can}_+$. Note that the centers of blowing up, being toroidal, 
automatically have normal crossings with $W^{can}_\pm \setmin U^{can}_\pm$. The
theorem follows. \qed

\section{Problems related to weak factorization}\label{Sec:problems}

\subsection{Strong factorization}\label{Sec:strong-factorization}
Despite our attempts, we have not  been able to use the methods of this
paper to prove the strong factorization conjecture, even assuming the toroidal
case holds true.

In the construction of the torific ideal in \ref{Sec:torific-construction}
 and the analysis of its blowing up
in \ref{Sec:torifying-property}  and \ref{Sec:lifting-toroidal}, the 
assumption of the cobordism $B_{a_i}$ being quasi-elementary is essential. 
One can extend the ideal over the entire cobordism $B$, for
instance  by
taking the Zariski closure of its zero scheme, but the behavior of this
extension (as well as others we have considered) along $B \setmin B_{a_i}$ is
problematic. 

The weak factorization theorem reduces the strong factorization conjecture to
the following problem:

\begin{problem} Let $X_1 \to X_2 \to\cdots \to X_n$ be a sequence of blowings
up with nonsingular centers, with $X_n$ nonsingular, and such that the center
of blowing up of $X_i \to X_{i+1}$ has normal crossings with the exceptional
divisor of $X_{i+1}\to X_n$. Let $Y \to X_n$ be a blowing up with nonsingular
center. Find a strong factorization of  the birational map $X_1 \das Y$. 
\end{problem}

We believe that at least the threefold case of this problem is tractable.

\subsection{Toroidalization}
\begin{problem}[Toroidalization] Let $\phi:X \rightarrow Y$
be a surjective proper morphism between complete nonsingular
 varieties over an algebraically closed field of characteristic 0.  Do there
 exist  
 sequences of blowings up with smooth centers $\nu_X:{\tilde X} \rightarrow X$
and $\nu_Y:{\tilde 
Y} \rightarrow Y$ so that the induced map ${\tilde \phi}:{\tilde X}
\das {\tilde Y}$ is a toroidal morphism?
Can such maps be chosen in a functorial manner, and in such a way that they
 preserve any open set where $\phi$ admits a toroidal structure?
\end{problem}

A similar conjecture was proposed in \cite{King1}. We note that the
toroidalization conjecture  concerns not only  
birational morphisms $\phi$ but also  generically finite morphisms or
  morphisms with $\dim X > \dim Y$.  The solution to the above
conjecture would reduce the strong factorization
conjecture to the toroidal case, simply by considering the case of a birational
morphism $\phi$ and 
then applying the toroidal case to ${\tilde \phi}$.  At
present the authors know  a complete proof only if either $\dim X =2$ (see
below), or $\dim Y=1$ (which follows immediately from resolution of
singularities, see \cite{KKMS}, II \S3).
Recently, S. D. Cutkosky announced a solution of the case $\dim X = 3,\, \dim Y
= 2$.  

The conjecture is false in positive charactersitics due to wild
ramifications. See, e.g., \cite{Cutkosky-Piltant}. 

One general result which we do know is the following.
\begin{theorem} Let $\phi:X \rightarrow Y$
be a proper surjective morphism between complete
 varieties over an algebraically closed field of characteristic 0. Then there
 exists a modification $\nu_X:{\tilde X} \rightarrow X$ 
 and a
 sequence of blowings up with smooth centers  $\nu_Y:{\tilde 
Y} \rightarrow Y$ so that the induced map ${\tilde \phi}:{\tilde X}
\das {\tilde Y}$ is a toroidal morphism.
\end{theorem}

{\bf Proof.} In \cite{Abramovich-Karu},  Theorem 2.1, it is shown that
modifications $\nu_X$ and $\nu_Y$ such that ${\tilde \phi}$ is toroidal exist,
assuming $X$ and $Y$ are projective and the generic fiber of $\phi$ is
geometrically integral. We can reduce to the projective case using Chow's
lemma. The case where the generic fiber is not geometrically integral is
resolved in the second author's thesis \cite{Karu-thesis}. Since the latter is
not widely available we give a similar argument here. The inductive proof of
\cite{Abramovich-Karu},  Theorem 2.1 reduces the problem to the case where
$\phi$ is generically finite. By Hironaka's flattening (or by taking a
resolution of the graph of $Y \das Hilb_Y(X)$), we may assume that $X \to Y$ is
finite.  Using resolution of
singularities, we may assume $Y$ is nonsingular and 
the branch locus is a normal crossings divisor. By normalizing $X$ we may
assume $X$ normal. Denoting the complement of the
branch locus by $U_Y$ and its inverse image in $X$ by $U_X$, Abhyankar's lemma
says that $U_X \subset X$ is a toroidal embedding and $X \to  Y $ is toroidal,
which is what we needed.

It remains to be shown that $\nu_Y$ can be chosen to be a sequence of blowings
up with nonsingular centers.  Let $Y \leftarrow Y' \to \tilde Y$ be a
resolution of indeterminacies of $Y \das \tilde Y$  and let $Y'' \to Y'$ be the
canonical principalization of the pullback of the ideal of the toroidal divisor
of $\tilde Y$. Let $X''\to Y'' \times_{\tilde Y}\tilde X$ be the normalization
of the dominant component. 
Then $Y'' \to Y$ is a sequence of   blowings
up with nonsingular centers. Applying \cite{Abramovich-Karu}, Lemma 6.2, we
see that $X'' \to Y''$ is still toroidal, which is what we needed.\qed

Since every proper birational morphism of nonsingular surfaces factors as a
sequence of point blowings up, we get:
\begin{corollary}
The toroidalization conjecture holds for a generically finite morphism $\phi: X
\to Y$ 
of surfaces.
\end{corollary}
In this case, it is not difficult to deduce that there exists a {\em minimal}
toroidalization (since the configuration of intermediate blowings up in $\tilde
X \to X$ or $\tilde Y \to Y$ forms a 
tree). This 
result has been proven in an algorithmic manner by Cutkosky and Piltant
\cite{Cutkosky-Piltant}. Similar statements can be found in
\cite{Akbulut-King}.

\end{document}